\documentclass{amsart}

\newtheorem{Theorem}{Theorem}[section]
\newtheorem{Proposition}[Theorem]{Proposition}
\newtheorem{Lemma}[Theorem]{Lemma}
\newtheorem{Corollary}[Theorem]{Corollary}
\theoremstyle{remark}
\newtheorem{Remark}[Theorem]{Remark}

\numberwithin{equation}{section}

\allowdisplaybreaks

\begin{document}

\title[Multilateral transformations of $q$-series]
{Multilateral transformations of $q$-series with quotients
of parameters that are nonnegative integer powers of $q$}
\author{Michael Schlosser}
\address{Department of Mathematics, The Ohio State University,
231 West 18th Avenue, Columbus, Ohio 43210, USA}
\email{mschloss@math.ohio-state.edu}
\urladdr{http://www.math.ohio-state.edu/\textasciitilde mschloss}
\date{February 21, 2001}

\begin{abstract}
We give multidimensional generalizations of several transformation
formulae for basic hypergeometric series of a specific type. Most of the
upper parameters of the series differ multiplicatively from corresponding
lower parameters by a nonnegative integer power of the base $q$. In one
dimension, formulae for such series have been found, in the $q\to 1$ case,
by B.~M.~Minton and P.~W.~Karlsson, and in the basic case by G.~Gasper, by
W.~C.~Chu, and more recently by the author. Our identities involve
multilateral basic hypergeometric series associated to the root system
$A_{r}$ (or equivalently, the unitary group $U\!(r+1)$).
\end{abstract}

\maketitle

\section{Introduction}

The theory of hypergeometric and basic hypergeometric (or $q$-hypergeometric)
series 
(cf. L.~J.~Slater~\cite{slater}, and G.~Gasper and M.~Rahman~\cite{grhyp})
contains numerous summation and transformation formulae.
Many of these appear in applications including number theory,
combinatorics, physics, representation theory, and computer algebra
(see e.g.\ G.~E.~Andrews~\cite{qandrews}).

One particular example is B.~M.~Minton's~\cite{minton} summation formula,
found in 1970, which is useful for simplifying sums that arise in certain
problems in theoretical physics (such as Racah coefficients). B.~M.~Minton's
formula is of special interest since it sums a specific hypergeometric
series with an abitrary number of parameters.
B.~M.~Minton derived his formula by expanding a hypergeometric series in terms
of other hypergeometric series, exploiting an identity already obtained by
C.~Fox~\cite{fox} in 1925. B.~M.~Minton iterated this expansion and suitably
specialized the parameters to successively evaluate the (inner) sums.
A condition on the parameters
of the specific hypergeometric series considered by B.~M.~Minton is that most
of the upper parameters differ from corresponding lower ones by a nonnegative
integer. B.~M.~Minton's result was slightly extended by
P.~W.~Karlsson~\cite{karlsson} who was using the same method.

In the early 1980's, G.~Gasper~\cite{gassum} found $q$-analogues of Karlsson
and Minton's results. In the basic case, the condition on the parameters
is that most of the upper parameters differ {\em multiplicatively} from
corresponding lower ones by a nonnegative integer power of $q$. 
G.~Gasper even extended his results to a transformation
formula~\cite[Eq.~(19)]{gassum}. For the above material,
see the exposition in G.~Gasper and M.~Rahman~\cite{grhyp}, in particular
Section~1.9, and Exercises~1.30 and 1.34.

Note that G.~Gasper and M.~Rahman~\cite{grhyp} use the terminology
``Karlsson--Minton" and ``$q$-Karlsson--Minton", respectively,
to denote the type of the series in question.
We are dropping this terminology in the present paper, since the work
is based on expanding a hypergeometric function in terms of another,
which has a longer history.
Instead, we introduce the acronyms IPD and $q$-IPD, respectively,
where IPD stands for ``{\bf I}ntegral {\bf P}arameter {\bf D}ifferences"
(motivated by the title of P.~W.~Karlsson's~\cite{karlsson} article),
see Section~\ref{secpre}.
It should be mentioned that expansions of hypergeometric series
in terms of other hypergeometric series have also been obtained
by J.~L.~Fields and J.~Wimp~\cite{fieldswimp}, by A.~Verma~\cite{averma},
and in more generality (concerning identities between general sequences),
by J.~L.~Fields and M.~E.~H.~Ismail~\cite{fieldsmail}.
Thus, as pointed out to us by Mourad Ismail~\cite{ismpriv},
one can easily write generalizations of the Karlsson--Minton
formulae to series involving partly hypergeometric
coefficients and partly general sequences.

By using an essentially different method, namely by partial fraction expansions,
W.~C.~Chu~\cite{chubs} generalized G.~Gasper's $q$-IPD type identities
further to a bilateral series transformation. 
In another article, G.~Gasper~\cite[Eq.~(5.13)]{gaseltf} found a new
summation for a very-well-poised basic hypergeometric series of $q$-IPD
type. Again, W.~C.~Chu~\cite{chuwp} extended G.~Gasper's result to a
summation for a very-well-poised {\em bilateral}
basic hypergeometric series.

Very recently, the author~\cite[Sec.~8]{schleltf} found even more general
identities of $q$-IPD type, by elementary manipulations of series,
using L.~J.~Slater's~\cite{slatertf} general transformations for
bilateral basic hypergeometric series. Already earlier 
J.~Haglund~\cite[pp.~415--416]{hagl} had discovered that
W.~C.~Chu's~\cite{chubs} bilateral transformation formula can be obtained by 
specializing L.~J.~Slater's~\cite{slatertf} general transformation for
${}_t\psi_t$ series.

In this article, we provide {\em multidimensional} extensions
of several specific transformation formulae of $q$-IPD type, in
particular, multivariate extensions of the identities in
Propositions~\ref{km3}, \ref{km2}, \ref{km1} and \ref{km0}.
These multivariate extensions involve multiple basic hypergeometric
series associated to the root system $A_{r-1}$ (or equivalently,
the unitary group $U\!(r)$). Such type of series are considered in the
work of R.~A.~Gustafson, S.~C.~Milne, and several other authors, see e.g.\
\cite{bhatmil}, \cite{geskratt} \cite{gusmult}, \cite{gus}, \cite{milmac},
\cite{miln1}, \cite{milne}, \cite{millil2}, \cite{milnew}, \cite{milschloss},
\cite{schlossmmi}, \cite{schlnammi}, and \cite{hypdet}.

As a matter of fact, there are unfortunately no suitable {\em multidimensional}
extensions of L.~J.~Slater's~\cite{slatertf} general transformation formulae
known (yet).
Thus, in higher dimensions we cannot specialize down from
such higher level identities. Instead we proceed from lower level
identities to systematically derive the upper level ones.
In this fashion, using certain $A_{r-1}$ summation theorems 
(from R.~A.~Gustafson~\cite{gusmult} and S.~C.~Milne~\cite{milne}),
elementary manipulation of series, and induction, we prove two 
multilateral transformations of $q$-IPD type,
namely Theorems~\ref{mchu} and \ref{mchua}. The first one of these,
Theorem~\ref{mchu}, involves {\em very-well-poised} multilateral series
(over $A_{r-1}$), and contains $r$-dimensional generalizations of 
W.~C.~Chu's~\cite[Theorem~2]{chuwp} and G.~Gasper's~\cite[Eq.~(5.13)]{gaseltf}
summations as special cases, see Corollaries~\ref{mchuc} and \ref{mgasc},
respectively. The other transformation formula
in Theorem~\ref{mchua}, involves multilateral series with an arbitrary
argument $z$. Four other multilateral transformations of
$q$-IPD type are derived by simpler means,
using tools developed in \cite{milschloss}, see Theorems~\ref{mkm1},
\ref{mkm1a}, \ref{mkm0}, and \ref{mkm0a}. 

In \cite[Theorem~6.4]{hypdet}, we already gave some multiple series
generalizations (associated to the root systems of classical type) of
W.~C.~Chu's~\cite{chubs} bilateral transformation.
The multiple series identities in \cite{hypdet}
were derived by using one-dimensional identities,
combined with certain determinant evaluations. In the same manner, one could
also deduce multilateral generalizations of L.~J.~Slater's~\cite{slatertf}
general transformation formulae, and in particular of the
$q$-IPD type transformations which were found in
\cite[Sec.~8]{schleltf}.
The identities one would obtain by this determinant method
would be not as deep as the ones derived in this article, though.

Our article is organized as follows.
In Section~\ref{secpre}, we introduce some standard notation for $q$-series
and basic hypergeometric series, and state several important one-dimensional
results. In Section~\ref{secmult}, we consider multiple series and
recollect some specific ingredients which we need in
Section~\ref{secmain} to state and prove our multilateral identities of
$q$-IPD type.

\section{Notation and one-dimensional results}\label{secpre}

In order to state and prove our theorems, we employ some standard $q$-series
notation (cf.\ G.~Gasper and M.~Rahman~\cite{grhyp}). For a complex number
$q$ with $0<|q|<1$, define the {\em $q$-shifted factorial} by
\begin{equation*}
(a;q)_{\infty}:=\prod_{j=0}^{\infty}(1-aq^j),
\end{equation*}
and
\begin{equation}\label{qpochdef}
(a;q)_k:=\frac{(a;q)_{\infty}}{(aq^k;q)_{\infty}},\qquad
\text{where $k$ is an integer}.
\end{equation}
Further, for brevity, we also employ the notation
\begin{equation*}
(a_1,\ldots,a_m;q)_k\equiv (a_1;q)_k\dots(a_m;q)_k,
\end{equation*}
where $k$ is an integer or infinity. Further, we utilize the notations
\begin{equation}\label{defhyp}
{}_t\phi_{t-1}\!\left[\begin{matrix}a_1,a_2,\dots,a_t\\
b_1,b_2,\dots,b_{t-1}\end{matrix}\,;q,z\right]:=
\sum _{k=0} ^{\infty}\frac {(a_1,a_2,\dots,a_t;q)_k}
{(q,b_1,\dots,b_{t-1};q)_k}z^k,
\end{equation}
and
\begin{equation}\label{defhypb}
{}_t\psi_t\!\left[\begin{matrix}a_1,a_2,\dots,a_t\\
b_1,b_2,\dots,b_t\end{matrix}\,;q,z\right]:=
\sum _{k=-\infty} ^{\infty}\frac {(a_1,a_2,\dots,a_t;q)_k}
{(b_1,b_2,\dots,b_t;q)_k}z^k,
\end{equation}
for {\em basic hypergeometric $_t\phi_{t-1}$ series}, and {\em bilateral basic
hypergeometric ${}_t\psi_t$ series}, respectively. Note that
G.~Gasper and M.~Rahman~\cite{grhyp} have more general definitions
for ${}_r\phi_s$ series and for ${}_r\psi_s$ series, but in this article
we are only really concerned with the case where $s=r-1$ for the
${}_r\phi_s$ series, and where $r=s$ for the ${}_r\psi_s$ series.

Clearly, a bilateral ${}_t\psi_t$ series becomes a
unilateral $_t\phi_{t-1}$ series if one of the lower parameters, say
$b_t$, is $q$ (or more generally, $q^j$ where $j$ is a positive integer).
This is because $(q;q)_k^{-1}=0$, for $k=-1, -2, \ldots$,
by definition \eqref{qpochdef}. In this case, the ${}_t\psi_t$ series
terminates naturally from below. On the other hand,
if in a  $_t\phi_{t-1}$ series one of the upper parameters, say $a_t$,
equals $q^{-n}$, where $n$ is a nonnegative integer, then the
$_t\phi_{t-1}$ series terminates naturally from above.
This is because $(q^{-n};q)_k=0$, for $k=n+1, n+2,\ldots$,
by definition \eqref{qpochdef}. Such a $_t\phi_{t-1}$ series terminates 
after $n+1$ terms.

The ratio test gives simple criteria of when the above series converge,
if they do not terminate. Remember that we assume $0<|q|<1$.
The $_t\phi_{t-1}$ series in \eqref{defhyp} converges absolutely in the radius
$|z|<1$, while the ${}_t\psi_t$ series in \eqref{defhypb} converges absolutely
in the annulus $|b_1\dots b_t/a_1\dots a_t|<|z|<1$. 

The classical theory of basic hypergeometric series consists of
several summation and transformation formulae involving $_t\phi_{t-1}$ series.
The classical summation theorems for terminating
$_3\phi_2$, $_6\phi_5$, and $_8\phi_7$ series require that the parameters
satisfy the additional condition of being either balanced and/or
very-well-poised. A $_t\phi_{t-1}$ basic hypergeometric series is called
{\em balanced} if $b_1\cdots b_{t-1}=a_1\cdots a_tq$ and $z=q$.
An $_t\phi_{t-1}$ series is {\em well-poised} if
$a_1q=a_2b_1=\cdots=a_tb_{t-1}$.
It is called {\em very-well-poised} if it is well-poised and if
$a_2=q\sqrt{a_1}$ and $a_3=-q\sqrt{a_1}$.
Note that  the factor
\begin{equation}\label{vwp}
\frac{(q\sqrt{a_1},-q\sqrt{a_1};q)_k}{(\sqrt{a_1},-\sqrt{a_1};q)_k}=
\frac {1-a_1q^{2k}}{1-a_1}
\end{equation}
appears in  a very-well-poised series.
The parameter $a_1$ is usually referred to as the
{\em special parameter} of such a series, and we call \eqref{vwp}
the {\em very-well-poised term} of the series.
Similarly, a bilateral $_t\psi_t$ basic
hypergeometric series is well-poised if
$a_1b_1=a_2b_2\cdots=a_tb_t$ and very-well-poised if, in addition,
$a_1=-a_2=qb_1=-qb_2$.

In our proofs in Section~\ref{secmain}, we often make use of some
elementary identities involving $q$-shifted factorials, listed
in G.~Gasper and M.~Rahman~\cite[Appendix~I]{grhyp}.

With the above notations for basic hypergeometric and bilateral basic
hypergeometric series, we are ready to
state some important (one-dimensional) summation formulae.

One of the most fundamental summation theorems in the theory of
(bilateral) basic hypergeometric series is W.~N.~Bailey's~\cite{bail66}
very-well-poised $_6\psi_6$ summation,
\begin{multline}\label{66gl}
{}_6\psi_6\!\left[\begin{matrix}q{\sqrt a}, -q{\sqrt a},b,c,d,e\\
{\sqrt a},-{\sqrt a},aq/b,aq/c,aq/d,aq/e\end{matrix}\,;
q,\frac {a^2q} {bcde}\right]\\
=\frac {(aq,aq/bc,aq/bd,aq/be,aq/cd,aq/ce,aq/de,q,q/a;q)_{\infty}}
{(aq/b,aq/c,aq/d,aq/e,q/b,q/c,q/d,q/e,a^2q/bcde;q)_{\infty}},
\end{multline}
provided the series either terminates, or $|q|<1$ and $|a^2q/bcde|<1$,
for convergence. For a simple proof of \eqref{66gl} using elementary
manipulations of series, see \cite{schlelsum}.

Another important summation is the terminating balanced
$q$-Pfaff--Saalsch\"utz summation (cf.~\cite[Eq.~(II.12)]{grhyp}),
\begin{equation}\label{32gl}
{}_3\phi_2\!\left[\begin{matrix}a,b,q^{-n}\\
c,abq^{1-n}/c\end{matrix}\,;
q,q\right]=
\frac{(c/a,c/b;q)_n}{(c,c/ab;q)_n}.
\end{equation}

S.~Ramanujan's $_1\psi_1$ summation (cf.~\cite[Eq.~(5.2.1)]{grhyp})
reads as follows,
\begin{equation}\label{11gl}
{}_1\psi_1\!\left[\begin{matrix}a\\b\end{matrix}\,;q,z\right]
=\frac{(q,b/a,az,q/az;q)_{\infty}}{(b,q/a,z,b/az;q)_{\infty}},
\end{equation}
provided the series either terminates, or $|q|<1$ and $|b/a|<|z|<1$,
for convergence.

Finally, the terminating $q$-binomial theorem is
(cf.~\cite[Eq.~(II.4)]{grhyp})
\begin{equation}\label{10tgl}
{}_1\phi_0\!\left[\begin{matrix}q^{-n}\\-\end{matrix}\,;q,z\right]
=(zq^{-n};q)_n.
\end{equation}
Note that \eqref{10tgl} is just the special case
$a\to q^{-n}$, $b\to q$ of \eqref{11gl}.

In this article, we prove multidimensional extensions (associated to the root
system $A_{r-1}$) of four transformations of $q$-IPD type,
namely Propositions~\ref{km3}, \ref{km2}, \ref{km1}, and \ref{km0}.
We need to explain our terminology first.

We say that a basic hypergeometric series is of
{\em $q$-IPD type} if there are $s$ upper parameters
$a_1,\dots,a_s$ and $s$ lower parameters $b_1,\dots,b_s$ such that
each $a_i$ differs from $b_i$ multiplicatively
by a nonnegative integer power of $q$,
i.e.\ $a_i=b_iq^{m_i}$, $m_i\ge 0$.
G.~Gasper~\cite{gassum} found some summation formulae for particular
basic hypergeometric series of such type. These were $q$-analogues
of formulae originally discovered by B.~M.~Minton~\cite{minton} and
P.~W.~Karlsson~\cite{karlsson}, using C.~Fox'~\cite{fox} expansion 
of a hypergeometric function in terms of other hypergeometric functions.
We call the series considered by B.~M.~Minton and P.~W.~Karlsson
to be of {\em IPD type}, where IPD stands for ``{\bf I}ntegral {\bf P}arameter
{\bf D}ifferences", motivated by the title of P.~W.~Karlsson's~\cite{karlsson}
article. G.~Gasper~\cite[Eq.~(19)]{gassum} also
extended his summations to a transformation formula. Later,
W.~C.~Chu~\cite{chubs} found {\em bilateral} summations and transformations of
$q$-IPD type, generalizing G.~Gasper's identities of \cite{gassum}.
In an expository paper, G.~Gasper~\cite[Eq.~(5.13)]{gaseltf}
derived a summation formula for a specific {\em very-well-poised}
basic hypergeometric series of $q$-IPD type.
His result was then generalized to a summation for bilateral series, again by
W.~C.~Chu~\cite[Theorem~2]{chuwp}. (It is maybe interesting that as application
W.~C.~Chu~\cite[Eq.~(5.25)]{chuwp} applied an inverse relation to his
bilateral summation and (re-)derived an important bibasic identity, actually
due to G.~Gasper and M.~Rahman~\cite[Eq.~(2.8)]{gasrahm}. This shows how
strongly seemingly different aspects in $q$-series are interconnected.)

In a recent article~\cite[Sec.~8]{schleltf}, the author found formulae of
$q$-IPD type covering all of the above $q$-IPD type identities as special cases.
In the following, we list the four transformation formulae from
\cite[Sec.~8]{schleltf} which we extend to higher dimensions.
The first one of these involves {\em very-well-poised}
bilateral basic hypergeometric series.

\begin{Proposition}[A bilateral very-well-poised $q$-IPD type transformation]
\label{km3}
Let $a$, $b$, $c$, $d$, $e$, $f$, and $h_1,\dots,h_s$
be indeterminate, let $m_1,\dots,m_s$
be nonnegative integers, let $|m|=\sum_{i=1}^s m_i$, and suppose
that the series in \eqref{km3gl} are well-defined. Then
\begin{multline}\label{km3gl}
{}_{6+2s}\psi_{6+2s}\!\Bigg[\begin{matrix}q\sqrt{a},-q\sqrt{a},b,c,d,e,\\
\sqrt{a},-\sqrt{a},\frac{a q}{b},\frac{a q}{c},\frac{a q}{d},\frac{a q}{e},
\end{matrix}\\
\begin{matrix}
h_1,\dots,h_s,\frac{aq^{1+m_1}}{h_1},\dots,\frac{aq^{1+m_s}}{h_s}\\
\frac{aq}{h_1},\dots,\frac{aq}{h_s},h_1q^{-m_1},\dots,h_sq^{-m_s}\end{matrix}
\,;q,\frac{a^{2}q^{1-|m|}}{bcde}\Bigg]\\
=\frac{(a,\frac q{a},\frac{fq}b,\frac{fq}c,\frac{fq}d,\frac{fq}e,
\frac{a q}{bf},\frac{a q}{cf},\frac{a q}{df},\frac{a q}{ef};q)_{\infty}}
{(\frac qb,\frac qc,\frac qd,\frac qe,\frac{a q}b,\frac{a q}c,
\frac{a q}d,\frac{a q}e,\frac{f^2q}{a},\frac{aq}{f^2};q)_{\infty}}
\prod_{i=1}^s\frac{(\frac{fq}{h_i},\frac{aq}{fh_i};q)_{m_i}}
{(\frac{aq}{h_i},\frac q{h_i};q)_{m_i}}\\\times
{}_{6+2s}\psi_{6+2s}\!\Bigg[\begin{matrix}\frac{qf}{\sqrt{a}},
-\frac{qf}{\sqrt{a}},\frac{bf}{a},\frac{cf}{a},\frac{df}{a},\frac{ef}{a},\\
\frac{f}{\sqrt{a}},-\frac{f}{\sqrt{a}},\frac{fq}{b},
\frac{fq}{c},\frac{fq}{d},\frac{fq}{e},
\end{matrix}\\
\begin{matrix}
\frac{fh_1}a,\dots,\frac{fh_s}a,\frac{fq^{1+m_1}}{h_1},\dots,
\frac{fq^{1+m_s}}{h_s}\\
\frac{fq}{h_1},\dots,\frac{fq}{h_s},\frac{fh_1q^{-m_1}}a,\dots,
\frac{fh_sq^{-m_s}}a\end{matrix}\,;q,
\frac{a^{2}q^{1-|m|}}{bcde}\Bigg],
\end{multline}
where the series either terminate, or $|a^{2}q^{1-|m|}/bcde|<1$,
for convergence.
\end{Proposition}
Note that $f$ does not appear on the left side of \eqref{km3gl}.

The special case $f\mapsto b$, $c\mapsto a/b$ of Proposition~\ref{km3}
is exactly W.~C.~Chu's summation in \cite[Theorem~2]{chuwp}. If we
specialize this summation then further by setting $e\mapsto a$ we arrive at
G.~Gasper's~\cite[Eq.~(5.13)]{gaseltf} summation.

The following transformation formula involves bilateral basic hypergeometric
series with an independent argument $z$.

\begin{Proposition}[A bilateral $q$-IPD type transformation]\label{km2}
Let $a$, $b$, $c$, $z$, and $h_1,\dots,h_s$ be indeterminate,
let $m_1,\dots,m_s$ be nonnegative integers, let $|m|=\sum_{i=1}^s m_i$,
and suppose that the series in \eqref{km2gl} are well-defined. Then
\begin{multline}\label{km2gl}
{}_{1+s}\psi_{1+s}\!\left[\begin{matrix}a,h_1q^{m_1},\dots,h_sq^{m_s}\\
b,h_1,\dots,h_s\end{matrix}\,;q,z\right]\\
=\frac{(c/a,bq/c,az,q/az;q)_{\infty}}
{(q/a,b,azq/c,c/az;q)_{\infty}}
\prod_{i=1}^s\frac{(h_iq/c;q)_{m_i}}{(h_i;q)_{m_i}}\\\times
{}_{1+s}\psi_{1+s}\!\left[\begin{matrix}aq/c,h_1q^{1+m_1}/c,\dots,
h_sq^{1+m_s}/c\\
bq/c,h_1q/c,\dots,h_sq/c\end{matrix}\,;q,z\right],
\end{multline}
where the series either terminate, or $|bq^{-|m|}/a|<|z|<1$,
for convergence.
\end{Proposition}
Note that $c$ does not appear on the left side of \eqref{km2gl}.

The next two transformations involve series whose argument depends
on the parameters.

\begin{Proposition}[A bilateral $q$-IPD type transformation]\label{km1}
Let $a$, $b$, $c$, $d$, $e$, and $h_1,\dots,h_s$ be indeterminate,
let $N$ be an arbitrary integer, $m_1,\dots,m_s$ be nonnegative integers,
let $|m|=\sum_{i=1}^s m_i$, and suppose that the series in
\eqref{km1gl} are well-defined. Then
\begin{multline}\label{km1gl}
{}_{2+s}\psi_{2+s}\!\left[\begin{matrix}a,b,h_1q^{m_1},\dots,h_sq^{m_s}\\
c,d,h_1,\dots,h_s\end{matrix}\,;q,\frac{eq^{-N}}{ab}\right]\\
=\left(\frac eq\right)^N\frac{(e/a,e/b,cq/e,dq/e;q)_{\infty}}
{(q/a,q/b,c,d;q)_{\infty}}
\prod_{i=1}^s\frac{(h_iq/e;q)_{m_i}}{(h_i;q)_{m_i}}\\\times
{}_{2+s}\psi_{2+s}\!\left[\begin{matrix}aq/e,bq/e,h_1q^{1+m_1}/e,\dots,
h_sq^{1+m_s}/e\\
cq/e,dq/e,h_1q/e,\dots,h_sq/e\end{matrix}\,;q,\frac{eq^{-N}}{ab}\right],
\end{multline}
where the series either terminate, or $|e/ab|<|q^N|<|eq^{|m|}/cd|$,
for convergence.
\end{Proposition}

If we reverse the ${}_{2+s}\psi_{2+s}$ series on the right side
of \eqref{km1gl}, we obtain
\begin{Proposition}[A bilateral $q$-IPD type transformation]\label{km0}
Let $a$, $b$, $c$, $d$, $e$, and $h_1,\dots,h_s$ be indeterminate,
let $N$ be an arbitrary integer, $m_1,\dots,m_s$ be nonnegative integers,
let $|m|=\sum_{i=1}^s m_i$, and suppose that the series in
\eqref{km0gl} are well-defined. Then
\begin{multline}\label{km0gl}
{}_{2+s}\psi_{2+s}\!\left[\begin{matrix}a,b,h_1q^{m_1},\dots,h_sq^{m_s}\\
c,d,h_1,\dots,h_s\end{matrix}\,;q,\frac{eq^{-N}}{ab}\right]\\
=\left(\frac eq\right)^N\frac{(e/a,e/b,cq/e,dq/e;q)_{\infty}}
{(q/a,q/b,c,d;q)_{\infty}}
\prod_{i=1}^s\frac{(h_iq/e;q)_{m_i}}{(h_i;q)_{m_i}}\\\times
{}_{2+s}\psi_{2+s}\!\left[\begin{matrix}e/c,e/d,e/h_1,\dots,e/h_s\\
e/a,e/b,eq^{-m_1}/h_1,\dots,eq^{-m_s}/h_s\end{matrix}\,;q,
\frac{cdq^{N-|m|}}e\right],
\end{multline}
where the series either terminate, or $|e/ab|<|q^N|<|eq^{|m|}/cd|$,
for convergence.
\end{Proposition}

The $e=aq$ case of Proposition~\ref{km0} reduces to
W.~C.~Chu's~\cite[Eq.~(15)]{chubs} transformation. If we specialize
the resulting transformation further by setting $c=q$ we obtain
G.~Gasper's~\cite[Eq.~(19)]{gassum} transformation.

Propositions~\ref{km3}, \ref{km2}, \ref{km1}, and \ref{km0} appeared as
Corollaries~8.6, 8.3, 8.2 and Equation~(8.8) in \cite{schleltf}.
They were originally derived as special cases from even more general
transformations for bilateral basic hypergeometric series of
$q$-IPD type.

\section{Preliminaries on multiple series}\label{secmult}

In general, we consider multiple series of the form
\begin{equation}\label{skgl}
\sum_{k_1,\dots,k_r=-\infty}^{\infty}S({\mathbf k}),
\end{equation}
where ${\mathbf k}=(k_1,\dots,k_r)$, which reduce to classical
(bilateral) basic hypergeometric series when $r=1$. We call
such a multiple basic hypergeometric series {\em balanced} if it reduces to 
a balanced series when $r=1$. Well-poised and very-well-poised
series are defined analogously\footnote{These definitions may seem far too
general but they are practical.}.
In case these series do not terminate from below, we also call such series
{\em multilateral} basic hypergeometric series.

In our particular cases, we also have
\begin{equation}\label{arvandy}
\prod_{1\le i<j<r} \left(\frac {z_iq^{k_i}-z_jq^{k_j}} {z_i-z_j}\right)
\end{equation} 
(or something similar), as a factor of $S({\mathbf k})$. A typical example
is the right side of \eqref{r66def}. Since we may associate \eqref{arvandy}
with the product side of the Weyl denominator formula for the
root system $A_{r-1}$ (see e.g.\ D.~Stanton~\cite{stan}), we call our series
$A_{r-1}$ basic hypergeometric series, in accordance with I.~M.~Gessel and
C.~Krattenthaler~\cite[Eq.~(7.1)]{geskratt}. Note that often in the literature
(e.g.\ \cite{bhat}, \cite{millil2}, \cite{milschloss}, \cite{schlossmmi},
\cite{schlnammi}) these $r$-dimensional series are (inprecisely) called
$A_r$ series instead of $A_{r-1}$ series.

For convenience, we frequently use the notation
$|{\mathbf k}|:= k_1+\dots+k_r$. Note that on the right side of
\eqref{r66def} we have (in addition to \eqref{arvandy})
\begin{equation}\label{mvwp}
\prod_{i=1}^r\left(\frac{1-az_iq^{k_i+|{\mathbf k}|}}{1-az_i}\right)
\end{equation}
appearing as a factor in the summand of the series. It is easy to see
that the $r=1$ case of \eqref{mvwp} essentially reduces to \eqref{vwp}.
To clarify the special appearance of the very-well-poised term
in the multidimensional case (and even in the one-dimensional) case,
it is useful to view the series in one higher dimension.
In particular, we can write
\begin{multline}\label{msvwp}
\prod_{1\le i<j<r} \left(\frac {z_iq^{k_i}-z_jq^{k_j}} {z_i-z_j}\right)
\prod_{i=1}^r\left(\frac{1-az_iq^{k_i+(k_1+\dots+k_r)}}{1-az_i}\right)\\
=q^{k_1+\dots+k_r}
\prod_{1\le i<j\le r+1} \left(\frac {z_iq^{k_i}-z_jq^{k_j}} {z_i-z_j}\right),
\end{multline}
where $z_{r+1}=1/a$ and $k_{r+1}=-(k_1+\dots+k_r)$.
Thus, some $A_{r-1}$ basic hypergeometric series identities are sometimes
better viewed as identites associated to the {\em affine} root system
$\tilde{A}_{r}$ (or, equivalently, the special unitary group $SU(r+1)$).
For such an example, see Remark~\ref{surem}.

Let $a$, $b_1,\dots,b_r$, $c$, $d$, $e_1,\dots,e_r$, $z_1,\dots,z_r$,
and $w$ be indeterminate.
For purpose of compact notation, we define for $r\ge 1$
\begin{multline}\label{r66def}
{}_6\Psi_6^{(r)}\!\left[a;b_1,\dots,b_r;c,d;e_1,\dots,e_r;z_1,\dots,z_r
\big|\,q,w\right]\\
:=\sum_{k_1,\dots,k_r=-\infty}^{\infty}
\Bigg(\prod_{1\le i<j\le r}
\left(\frac{z_iq^{k_i}-z_jq^{k_j}}{z_i-z_j}\right)
\prod_{i=1}^r\left(\frac{1-az_iq^{k_i+|{\mathbf k}|}}{1-az_i}\right)\\\times
\prod_{i,j=1}^r\frac{(b_jz_i/z_j;q)_{k_i}}{(az_iq/e_jz_j;q)_{k_i}}
\prod_{i=1}^r\frac{(e_iz_i;q)_{|{\mathbf k}|}}{(az_iq/b_i;q)_{|{\mathbf k}|}}
\\\times
\prod_{i=1}^r\frac{(cz_i;q)_{k_i}}{(az_iq/d;q)_{k_i}}\cdot
\frac{(d;q)_{|{\mathbf k}|}}{(aq/c;q)_{|{\mathbf k}|}}\,
w^{|{\mathbf k}|}\Bigg).
\end{multline}

The above $_6\Psi_6^{(r)}$ series is an $r$-dimensional $_6\psi_6$ series
(which reduces to a classical very-well-posied $_6\psi_6$ when $r=1$).

For convenience, we sometimes use capital letters to abbreviate the ($r$-fold)
products of certain variables. Specifically, in this article we use
$A\equiv a_1\cdots a_r$, $B\equiv b_1\cdots b_r$, $C\equiv c_1\cdots c_r$,
$E\equiv e_1\cdots e_r$, and $F\equiv f_1\cdots f_r$, respectively. 

In our derivation of the multilateral $q$-IPD type
transformation in Theorem~\ref{mchu}
we utilize the following $r$-dimensional generalization of
W.~N.~Bailey's summation formula in \eqref{66gl}.
\begin{Theorem}[(Gustafson) An $A_{r-1}$ $_6\psi_6$ summation]\label{r66}
Let $a$, $b_1,\dots,b_r$, $c$, $d$, $e_1,\dots,e_r$, and $z_1,\dots,z_r$
be indeterminate, let $r\ge 1$, and suppose that
none of the denominators in \eqref{r66gl} vanishes. Then
\begin{multline}\label{r66gl}
{}_6\Psi_6^{(r)}\!\left[a;b_1,\dots,b_r;c,d;e_1,\dots,e_r;z_1,\dots,z_r
\Big|\,q,\frac{a^{r+1}q}{BcdE}\right]\\
=\frac{(aq/Bc,a^rq/dE,aq/cd;q)_{\infty}}
{(a^{r+1}q/BcdE,aq/c,q/d;q)_{\infty}}
\prod_{i,j=1}^r\frac{(az_iq/b_ie_jz_j,qz_i/z_j;q)_{\infty}}
{(qz_i/b_iz_j,az_iq/e_jz_j;q)_{\infty}}\\\times
\prod_{i=1}^r\frac{(aq/ce_iz_i,az_iq/b_id,az_iq,q/az_i;q)_{\infty}}
{(az_iq/b_i,q/e_iz_i,q/cz_i,az_iq/d;q)_{\infty}},
\end{multline}
provided $|a^{r+1}q/BcdE|<1$.
\end{Theorem}

\begin{Remark}\label{surem}
Using \eqref{msvwp}, the multilateral identity in \eqref{r66gl} can also be
written in a more compact form. We then have
R.~A.~Gustafson's~\cite[Theorem~1.15]{gusmult} $\tilde{A}_{r}$ $_6\psi_6$
summation:
Let $a_1,\dots,a_{r+1}$, $b_1,\dots,b_{r+1}$, and $z_1,\dots,z_{r+1}$
be indeterminate, let $r\ge 1$, and suppose that none of the
denominators in \eqref{an1psi1cglN} vanishes. Then
\begin{multline}\label{an1psi1cglN}
\sum_{\begin{smallmatrix}-\infty\le k_1,\dots,k_{r+1}\le\infty\\
k_1+\dots+k_{r+1}=0\end{smallmatrix}}
\prod_{1\le i<j\le r+1}\left(\frac {z_iq^{k_i}-z_jq^{k_j}}{z_i-z_j}\right)
\prod_{i,j=1}^{r+1}\frac{(a_jz_i/z_j;q)_{k_i}}{(b_jz_i/z_j;q)_{k_i}}\\
=\frac{(b_1\dots b_{r+1}q^{-r},q/a_1\dots a_{r+1};q)_{\infty}}
{(q,b_1\dots b_{r+1}q^{-r}/a_1\dots a_{r+1};q)_{\infty}}
\prod_{i,j=1}^{r+1}
\frac{(qz_i/z_j,b_jz_i/a_iz_j;q)_{\infty}}
{(b_jz_i/z_j,z_iq/a_iz_j;q)_{\infty}},
\end{multline}
provided $|b_1\dots b_{r+1}q^{-r}/a_1\dots a_{r+1}|<1$.
It is not difficult to see that \eqref{an1psi1cglN} and \eqref{r66gl}
are equivalent.
\end{Remark}

We also need the following $r$-dimensional
generalization of the terminating $q$-Pfaff--Saalsch\"utz
summation from S.~C.~Milne~\cite[Theorem~4.15]{milne}.
\begin{Theorem}[(Milne) An $A_{r-1}$ terminating $_3\phi_2$
summation]\label{r32}
Let $a_1,\dots,a_r$, $b$, $c$, and $x_1,\dots,x_r$, be indeterminate, let
$N$ be a nonnegative integer, let $r\ge 1$, and suppose that
none of the denominators in \eqref{r32gl} vanishes. Then
\begin{multline}\label{r32gl}
\sum_{\begin{smallmatrix}k_1,\dots,k_r\ge 0\\
0\le |{\mathbf k}|\le N\end{smallmatrix}}
\Bigg(\prod_{1\le i<j\le r}
\left(\frac{x_iq^{k_i}-x_jq^{k_j}}{x_i-x_j}\right)
\prod_{i,j=1}^r\frac{(a_jx_i/x_j;q)_{k_i}}{(qx_i/x_j;q)_{k_i}}\\\times
\prod_{i=1}^r\frac{(bx_i;q)_{k_i}}{(cx_i;q)_{k_i}}\cdot
\frac{(q^{-N};q)_{|{\mathbf k}|}}{(a_1\dots a_rbq^{1-N}/c;q)_{|{\mathbf k}|}}
\,q^{|{\mathbf k}|}\Bigg)\\
=\frac{(c/b;q)_N}{(c/a_1\dots a_rb;q)_N}
\prod_{i=1}^r\frac{(cx_i/a_i;q)_N}{(cx_i;q)_N}.
\end{multline}
\end{Theorem}

The $r=1$ case of \eqref{r32gl} clearly reduces to \eqref{32gl}.

In our derivation of the multilateral $q$-IPD type
transformation in Theorem~\ref{mchua}
we utilize R.~A.~Gustafson's~\cite[Theorem~1.17]{gusmult} $A_{r-1}$
extension of S.~Ramanujan's $_1\psi_1$ summation~\eqref{11gl}.

\begin{Theorem}[(Gustafson) An $A_{r-1}$ $_1\psi_1$ summation]\label{r11}
Let $a_1,\dots,a_r$, $b_1,\dots,$ $b_r$, $x_1,\dots,x_r$, and $z$ be
indeterminate, let $r\ge 1$, and
suppose that none of the denominators in \eqref{r11gl} vanishes. Then
\begin{multline}\label{r11gl}
\sum_{k_1,\dots,k_r=-\infty}^{\infty}
\prod_{1\le i<j\le r}\left(\frac{x_iq^{k_i}-x_jq^{k_j}}{x_i-x_j}\right)
\prod_{i,j=1}^r\frac{(a_jx_i/x_j;q)_{k_i}}{(b_jx_i/x_j;q)_{k_i}}\,
z^{|{\mathbf k}|}\\
=\frac{(Az,q/Az;q)_{\infty}}{(z,Bq^{1-r}/Az;q)_{\infty}}
\prod_{i,j=1}^r\frac{(b_jx_i/a_ix_j,qx_i/x_j;q)_{\infty}}
{(qx_i/a_ix_j,b_jx_i/x_j;q)_{\infty}},
\end{multline}
where $|Bq^{1-r}/A|<|z|<1$.
\end{Theorem}

Further, we make use of the following terminating $q$-binomial theorem
from S.~C.~Milne~\cite[Theorem~5.46]{milne},
which is a multiple extension of \eqref{10tgl}.
\begin{Theorem}[(Milne) An $A_{r-1}$ terminating $q$-binomial
theorem]\label{r10}
Let $x_1,\dots,$ $x_r$, and $z$ be indeterminate, let $n_1,\dots,n_r$
be nonnegative integers, let $r\ge 1$, and
suppose that none of the denominators in \eqref{r10gl} vanishes. Then
\begin{multline}\label{r10gl}
\sum_{\begin{smallmatrix}0\le k_i\le n_i\\
i=1,\dots,r\end{smallmatrix}}
\Bigg(\prod_{1\le i<j\le r}\left(\frac{x_iq^{k_i}-x_jq^{k_j}}{x_i-x_j}\right)
\prod_{i,j=1}^r\frac{(q^{-n_j}x_i/x_j;q)_{k_i}}
{(qx_i/x_j;q)_{k_i}}\prod_{i=1}^rx_i^{k_i}\\\times
q^{-{\binom{|{\mathbf k}|}2}+\sum_{i=1}^r{\binom{k_i}2}}z^{|{\mathbf k}|}\Bigg)
=\prod_{i=1}^r(zx_iq^{-|{\mathbf n}|};q)_{n_i}.
\end{multline}
\end{Theorem}

In Section~\ref{secmain}, we also give two multiple series extensions each
of Propositions~\ref{km1} and \ref{km0}, see
Theorems~\ref{mkm1}, \ref{mkm1a}, \ref{mkm0}, and \ref{mkm0a}.
These $A_{r-1}$ extensions are not as deep as those in Theorems~\ref{mchu}
or \ref{mchua}. In our derivations, we make use of Lemmas~4.3 and 4.9 from
\cite{milschloss}, displayed as follows:

\begin{Lemma}\label{lem1}
Let $b_1,\dots,b_r$ and $x_1,\dots,x_r$ be indeterminate, let $r\ge 1$,
and suppose that none of the denominators in \eqref{lem1gl} vanishes.
Then, if $f(n)$ is an arbitrary function of integers $n$, we have
\begin{multline}\label{lem1gl}
\sum_{n=-\infty}^{\infty}\frac{f(n)}{(Bq^{1-r};q)_n}
=\frac{(q;q)_{\infty}}{(Bq^{1-r};q)_{\infty}}
\prod_{i,j=1}^r\frac{(b_jx_i/x_j;q)_{\infty}}
{(qx_i/x_j;q)_{\infty}}\\\times
\sum_{k_1,\dots,k_r=-\infty}^{\infty}
\Bigg(\prod_{1\le i<j\le r}\left(\frac {x_iq^{k_i}-x_jq^{k_j}}{x_i-x_j}\right)
\prod_{i,j=1}^r(b_jx_i/x_j;q)_{k_i}^{-1}
\prod_{i=1}^rx_i^{rk_i-|{\mathbf k}|}\\\times
(-1)^{(r-1)|{\mathbf k}|}
q^{-\binom{|{\mathbf k}|}2+r\sum_{i=1}^r\binom{k_i}2}
\cdot f(|{\mathbf k}|)\Bigg),
\end{multline}
provided the series converge.
\end{Lemma}

\begin{Lemma}\label{lem2}
Let $a_1,\dots,a_r$, $b_1,\dots,b_r$, and $x_1,\dots,x_r$ be indeterminate,
let $r\ge 1$, and suppose that none of the denominators in
\eqref{lem2gl} vanishes.
Then, if $g(n)$ is an arbitrary function of integers $n$, we have
\begin{multline}\label{lem2gl}
\sum_{n=-\infty}^{\infty}\frac{(A;q)_n}{(Bq^{1-r};q)_n}\,g(n)\\
=\frac{(q,Bq^{1-r}/A;q)_{\infty}}
{(Bq^{1-r},q/A;q)_{\infty}}
\prod_{i,j=1}^r\frac{(b_jx_i/x_j,
x_iq/a_ix_j;q)_{\infty}}
{(qx_i/x_j,b_jx_i/a_ix_j;q)_{\infty}}\\\times
\sum_{k_1,\dots,k_r=-\infty}^{\infty}
\prod_{1\le i<j\le r}\left(\frac {x_iq^{k_i}-x_jq^{k_j}}{x_i-x_j}\right)
\prod_{i,j=1}^r\frac{(a_jx_i/x_j;q)_{k_i}}
{(b_jx_i/x_j;q)_{k_i}}
\cdot g(|{\mathbf k}|),
\end{multline}
provided the series converge.
\end{Lemma}

\section{Multilateral identities of $q$-IPD type}
\label{secmain}

Here we give six new multilateral transformations of
$q$-IPD type, extending the $q$-IPD type transformations of
Propositions~\ref{km3}, \ref{km2}, \ref{km1}, and \ref{km0} to
higher dimensions. The transformation formula in
Theorem~\ref{mchu}, which generalizes Proposition~\ref{km3},
involves multiple series very-well-poised over the root system $A_{r-1}$.
A special case of that theorem is given as Corollary~\ref{mchuc}, which is a
multilateral summation formula extending
W.~C.~Chu's~\cite[Theorem~2]{chuwp} bilateral summation to $r$-dimensions.
A further specialization gives a multiple extension of
G.~Gasper's~\cite[Eq.~(5.13)]{gaseltf} very-well-poised summation,
see Corollary~\ref{mgasc}.
In Theorem~\ref{mchua} we provide an $A_{r-1}$ extension of
Proposition~\ref{km2}. The interesting feature about that transformation
is that it involves multilateral series with an {\em independent}
argument $z$ (subject to convergence), similar to
the case of S.~Ramanujan's $_1\psi_1$ summation~\eqref{11gl}
and its extension in Theorem~\ref{r11}.
We were, unfortunately, not able to give multidimensional extensions of
Propositions~\ref{km1} or \ref{km0} which are as deep as
Theorems~\ref{mchu} and \ref{mchua}. Instead, we derive multiple extensions
of a {\em simpler} type, using Lemmas~\ref{lem1} and \ref{lem2}.
Theorems~\ref{mkm1} and \ref{mkm1a} are simple $A_{r-1}$ extensions
of Proposition~\ref{km1}, while Theorems~\ref{mkm0} and \ref{mkm0a} are
simple $A_{r-1}$ extensions of Proposition~\ref{km0}.
Of course, by the same method one could also derive simple multilateral
generalizations of the $q$-IPD type
transformations in Propositions~\ref{km3} and \ref{km2}.
However, we decided to derive the identities of simpler type only
in the cases were we were unable to find corresponding deeper ones.

For our derivation of Theorem~\ref{mchu}, we need the following lemma,
which is easily established by applying Theorem~\ref{r66} twice.
\begin{Lemma}\label{66lem}
Let $a$, $b_1,\dots,b_r$, $c$, $d$, $e_1,\dots,e_r$, $f_1,\dots,f_r$,
and $z_1,\dots,z_r$ be indeterminate, let $r\ge 1$, and suppose that
none of the denominators in \eqref{66lemgl} vanishes. Then
\begin{multline}\label{66lemgl}
{}_6\Psi_6^{(r)}\!\left[a;b_1,\dots,b_r;c,d;e_1,\dots,e_r;z_1,\dots,z_r
\Big|\,q,\frac{a^{r+1}q}{BcdE}\right]\\
=\frac{(Fq/d,aq/cF;q)_{\infty}}{(aq/c,q/d;q)_{\infty}}
\prod_{i,j=1}^r\frac{(qz_i/z_j,az_iq/e_jf_iz_j,f_jz_iq/b_iz_j;q)_{\infty}}
{(qf_jz_i/f_iz_j,qz_i/b_iz_j,az_iq/e_jz_j;q)_{\infty}}\\\times
\prod_{i=1}^r
\frac{(az_iq,q/az_i,Fq/e_iz_i,az_iq/b_iF,az_iq/df_i,f_iq/cz_i;q)_{\infty}}
{(Ff_iq/az_i,az_iq/Ff_i,az_iq/b_i,q/e_iz_i,q/cz_i,az_iq/d;q)_{\infty}}\\\times
{}_6\Psi_6^{(r)}\!\left[\frac Fa;\frac{e_1f_1}a,\dots,\frac{e_rf_r}a;
\frac da,\frac{cF}a;\frac{b_1F}{af_1},\dots,\frac{b_rF}{af_r};
\frac{f_1}{z_1},\dots,\frac{f_r}{z_r}
\bigg|\,q,\frac{a^{r+1}q}{BcdE}\right],
\end{multline}
provided $|a^{r+1}q/BcdE|<1$.
\end{Lemma}

Now, for compact notation, let us extend definition \eqref{r66def} by
introducing additional indeterminates
$g_1,\dots,g_s$ and $h_1,\dots,h_s$:
\begin{multline}
{}_{6+2s}\Psi_{6+2s}^{(r)}\big[a;b_1,\dots,b_r;c,d;e_1,\dots,e_r;
z_1,\dots,z_r;\\
\left. g_1,\dots,g_s;h_1,\dots,h_s\big|\,q,w\right]\\
:=\sum_{k_1,\dots,k_r=-\infty}^{\infty}
\Bigg(\prod_{1\le i<j\le r}
\left(\frac{z_iq^{k_i}-z_jq^{k_j}}{z_i-z_j}\right)
\prod_{i=1}^r\left(\frac{1-az_iq^{k_i+|{\mathbf k}|}}{1-az_i}\right)\\\times
\prod_{i,j=1}^r\frac{(b_jz_i/z_j;q)_{k_i}}{(az_iq/e_jz_j;q)_{k_i}}
\prod_{i=1}^r\frac{(e_iz_i;q)_{|{\mathbf k}|}}{(az_iq/b_i;q)_{|{\mathbf k}|}}
\\\times
\prod_{i=1}^r\frac{(cz_i,g_1z_i,\dots,g_sz_i;q)_{k_i}}
{(az_iq/d,az_iq/h_1,\dots,az_iq/h_s;q)_{k_i}}\\\times
\frac{(d,h_1,\dots,h_s;q)_{|{\mathbf k}|}}
{(aq/c,aq/g_1,\dots,aq/g_s;q)_{|{\mathbf k}|}}\,
w^{|{\mathbf k}|}\Bigg).
\end{multline}

We have
\begin{Theorem}[A multilateral very-well-poised $A_{r-1}$ $q$-IPD type
transformation]\label{mchu}
Let $a$, $b_1,\dots,b_r$, $c$, $d$, $e_1,\dots,e_r$, $f_1,\dots,f_r$,
$z_1,\dots,z_r$, and $h_1,\dots,h_s$ be indeterminate, let
$N_1,\dots,N_s$ be nonnegative integers, let $|N|=\sum_{i=1}^s N_i$,
$r\ge 1$, and
suppose that none of the denominators in \eqref{mchugl} vanishes. Then
\begin{multline}\label{mchugl}
{}_{6+2s}\Psi_{6+2s}^{(r)}\bigg[a;b_1,\dots,b_r;c,d;e_1,\dots,e_r;
z_1,\dots,z_r;\\
\left.\frac{aq^{1+N_1}}{h_1},\dots,\frac{aq^{1+N_s}}{h_s};h_1,\dots,h_s
\bigg|\,q,\frac{a^{r+1}q^{1-|N|}}{BcdE}\right]\\
=\prod_{j=1}^s\left[\frac{(Fq/h_j;q)_{N_j}}{(q/h_j;q)_{N_j}}
\prod_{i=1}^r\frac{(az_iq/f_ih_j;q)_{N_j}}{(az_iq/h_j;q)_{N_j}}\right]\\\times
\frac{(Fq/d,aq/cF;q)_{\infty}}{(aq/c,q/d;q)_{\infty}}
\prod_{i,j=1}^r\frac{(qz_i/z_j,az_iq/e_jf_iz_j,f_jz_iq/b_iz_j;q)_{\infty}}
{(qf_jz_i/f_iz_j,qz_i/b_iz_j,az_iq/e_jz_j;q)_{\infty}}\\\times
\prod_{i=1}^r
\frac{(az_iq,q/az_i,Fq/e_iz_i,az_iq/b_iF,az_iq/df_i,f_iq/cz_i;q)_{\infty}}
{(Ff_iq/az_i,az_iq/Ff_i,az_iq/b_i,q/e_iz_i,q/cz_i,az_iq/d;q)_{\infty}}\\\times
{}_{6+2s}\Psi_{6+2s}^{(r)}\bigg[\frac Fa;\frac{e_1f_1}a,\dots,\frac{e_rf_r}a;
\frac da,\frac{cF}a;\frac{b_1F}{af_1},\dots,\frac{b_rF}{af_r};
\frac{f_1}{z_1},\dots,\frac{f_r}{z_r};\\
\left.\frac{h_1}a,\dots,\frac{h_s}a;
\frac{Fq^{1+N_1}}{h_1},\dots,\frac{Fq^{1+N_s}}{h_s}
\bigg|\,q,\frac{a^{r+1}q^{1-|N|}}{BcdE}\right],
\end{multline}
provided $|a^{r+1}q^{1-|N|}/BcdE|<1$.
\end{Theorem}

\begin{proof}
We proceed by induction on $s$. For $s=0$ \eqref{mchugl} is true
by Lemma~\ref{66lem}. So, suppose that the transformation is already shown
for $s \mapsto s-1$. Then, by using some elementary identities from
\cite[Appendix~I]{grhyp},
\begin{multline*}
{}_{6+2s}\Psi_{6+2s}^{(r)}\bigg[a;b_1,\dots,b_r;c,d;e_1,\dots,e_r;
z_1,\dots,z_r;\\
\left.\frac{aq^{1+N_1}}{h_1},\dots,\frac{aq^{1+N_s}}{h_s};h_1,\dots,h_s
\bigg|\,q,\frac{a^{r+1}q^{1-|N|}}{BcdE}\right]\\
=\sum_{k_1,\dots,k_r=-\infty}^{\infty}
\Bigg(\prod_{1\le i<j\le r}
\left(\frac{z_iq^{k_i}-z_jq^{k_j}}{z_i-z_j}\right)
\prod_{i=1}^r\left(\frac{1-az_iq^{k_i+|{\mathbf k}|}}{1-az_i}\right)\\\times
\prod_{i,j=1}^r\frac{(b_jz_i/z_j;q)_{k_i}}{(az_iq/e_jz_j;q)_{k_i}}
\prod_{i=1}^r\frac{(e_iz_i;q)_{|{\mathbf k}|}}{(az_iq/b_i;q)_{|{\mathbf k}|}}
\\\times\prod_{i=1}^r
\frac{(cz_i,az_iq^{1+N_1}/h_1,\dots,az_iq^{1+N_{s-1}}/h_{s-1};q)_{k_i}}
{(az_iq/d,az_iq/h_1,\dots,az_iq/h_{s-1};q)_{k_i}}\\\times
\frac{(d,h_1,\dots,h_{s-1};q)_{|{\mathbf k}|}}
{(aq/c,h_1q^{-N_1},\dots,h_{s-1}q^{-N_{s-1}};q)_{|{\mathbf k}|}}
\left(\frac{a^{r+1}q^{1-|N|}}{BcdE}\right)^{|{\mathbf k}|}\\
\times\frac{(h_s;q)_{|{\mathbf k}|}}
{(h_sq^{-N_s};q)_{|{\mathbf k}|}}
\prod_{i=1}^r\frac{(az_iq^{1+N_s}/h_s;q)_{k_i}}
{(az_iq/h_s;q)_{k_i}}\Bigg)\\
=\frac{(a^rq/Eh_s;q)_{N_s}}{(q/h_s;q)_{N_s}}
\prod_{i=1}^r\frac{(e_iz_iq/h_s;q)_{N_s}}{(az_iq/h_s;q)_{N_s}}\\\times
\sum_{k_1,\dots,k_r=-\infty}^{\infty}
\Bigg(\prod_{1\le i<j\le r}
\left(\frac{z_iq^{k_i}-z_jq^{k_j}}{z_i-z_j}\right)
\prod_{i=1}^r\left(\frac{1-az_iq^{k_i+|{\mathbf k}|}}{1-az_i}\right)\\\times
\prod_{i,j=1}^r\frac{(b_jz_i/z_j;q)_{k_i}}{(az_iq/e_jz_j;q)_{k_i}}
\prod_{i=1}^r\frac{(e_iz_i;q)_{|{\mathbf k}|}}{(az_iq/b_i;q)_{|{\mathbf k}|}}
\\\times\prod_{i=1}^r
\frac{(cz_i,az_iq^{1+N_1}/h_1,\dots,az_iq^{1+N_{s-1}}/h_{s-1};q)_{k_i}}
{(az_iq/d,az_iq/h_1,\dots,az_iq/h_{s-1};q)_{k_i}}\\\times
\frac{(d,h_1,\dots,h_{s-1};q)_{|{\mathbf k}|}}
{(aq/c,h_1q^{-N_1},\dots,h_{s-1}q^{-N_{s-1}};q)_{|{\mathbf k}|}}
\left(\frac{a^{r+1}q^{1-(N_1+\dots +N_{s-1})}}{BcdE}\right)^{|{\mathbf k}|}\\
\times\frac{(q^{1-|{\mathbf k}|}/h_s;q)_{N_s}}
{(a^rq/Eh_s;q)_{N_s}}
\prod_{i=1}^r\frac{(az_iq^{1+k_i}/h_s;q)_{N_s}}
{(e_iz_iq/h_s;q)_{N_s}}\Bigg).
\end{multline*}
Now we expand the last factors (those involving $(\cdot;q)_{N_s}$) by applying
the $a_i\mapsto e_iq^{-k_i}/a$, $b\mapsto q^{|{\mathbf k}|}$, $c\mapsto q/h_s$,
$x_i\mapsto e_iz_i$, $i=1,\dots,r$, and $N\mapsto N_s$
case of the $A_{r-1}$ $q$-Pfaff--Saalsch\"utz summation in Theorem~\ref{r32}.
We obtain
\begin{multline*}
\frac{(a^rq/Eh_s;q)_{N_s}}{(q/h_s;q)_{N_s}}
\prod_{i=1}^r\frac{(e_iz_iq/h_s;q)_{N_s}}{(az_iq/h_s;q)_{N_s}}\\\times
\sum_{k_1,\dots,k_r=-\infty}^{\infty}
\Bigg(\prod_{1\le i<j\le r}
\left(\frac{z_iq^{k_i}-z_jq^{k_j}}{z_i-z_j}\right)
\prod_{i=1}^r\left(\frac{1-az_iq^{k_i+|{\mathbf k}|}}{1-az_i}\right)\\\times
\prod_{i,j=1}^r\frac{(b_jz_i/z_j;q)_{k_i}}{(az_iq/e_jz_j;q)_{k_i}}
\prod_{i=1}^r\frac{(e_iz_i;q)_{|{\mathbf k}|}}{(az_iq/b_i;q)_{|{\mathbf k}|}}
\\\times\prod_{i=1}^r
\frac{(cz_i,az_iq^{1+N_1}/h_1,\dots,az_iq^{1+N_{s-1}}/h_{s-1};q)_{k_i}}
{(az_iq/d,az_iq/h_1,\dots,az_iq/h_{s-1};q)_{k_i}}\\\times
\frac{(d,h_1,\dots,h_{s-1};q)_{|{\mathbf k}|}}
{(aq/c,h_1q^{-N_1},\dots,h_{s-1}q^{-N_{s-1}};q)_{|{\mathbf k}|}}
\left(\frac{a^{r+1}q^{1-(N_1+\dots +N_{s-1})}}{BcdE}\right)^{|{\mathbf k}|}\\
\times\sum_{\begin{smallmatrix}l_1,\dots,l_r\ge 0\\
0\le |{\mathbf l}|\le N_s\end{smallmatrix}}
\prod_{1\le i<j\le r}
\left(\frac{e_iz_iq^{l_i}-e_jz_jq^{l_j}}{e_iz_i-e_jz_j}\right)
\prod_{i,j=1}^r\frac{(e_iz_iq^{-k_j}/az_j;q)_{l_i}}
{(qe_iz_i/e_jz_j;q)_{l_i}}\\\times
\prod_{i=1}^r\frac{(ez_iq^{|{\mathbf k}|};q)_{l_i}}{(e_iz_iq/h_s;q)_{l_i}}
\cdot\frac{(q^{-N_s};q)_{|{\mathbf l}|}}{(Eh_sq^{-N_s}/a^r;q)_{|{\mathbf l}|}}
\,q^{|{\mathbf l}|}\Bigg)\\\
=\frac{(a^rq/Eh_s;q)_{N_s}}{(q/h_s;q)_{N_s}}
\prod_{i=1}^r\frac{(e_iz_iq/h_s;q)_{N_s}}{(az_iq/h_s;q)_{N_s}}\\\times
\sum_{\begin{smallmatrix}l_1,\dots,l_r\ge 0\\
0\le |{\mathbf l}|\le N_s\end{smallmatrix}}
\Bigg(\prod_{1\le i<j\le r}
\left(\frac{e_iz_iq^{l_i}-e_jz_jq^{l_j}}{e_iz_i-e_jz_j}\right)
\prod_{i,j=1}^r\frac{(e_iz_i/az_j;q)_{l_i}}
{(qe_iz_i/e_jz_j;q)_{l_i}}\\\times
\prod_{i=1}^r\frac{(ez_i;q)_{l_i}}{(e_iz_iq/h_s;q)_{l_i}}
\cdot\frac{(q^{-N_s};q)_{|{\mathbf l}|}}{(Eh_sq^{-N_s}/a^r;q)_{|{\mathbf l}|}}
\,q^{|{\mathbf l}|}\\\times
{}_{6+2(s-1)}\Psi_{6+2(s-1)}^{(r)}\bigg[a;b_1,\dots,b_r;c,d;
e_1q^{l_1},\dots,e_rq^{l_r};
z_1,\dots,z_r;\\
\frac{aq^{1+N_1}}{h_1},\dots,\frac{aq^{1+N_{s-1}}}{h_{s-1}};
h_1,\dots,h_{s-1}\bigg|\,q,
\frac{a^{r+1}q^{1-(N_1+\dots+N_{s-1})-|{\mathbf l}|}}{BcdE}\bigg]\Bigg).
\end{multline*}
By the $e_i\mapsto e_iq^{l_i}$, $i=1,\dots,r$, case of the inductive
hypothesis we obtain
\begin{multline*}
\frac{(a^rq/Eh_s;q)_{N_s}}{(q/h_s;q)_{N_s}}
\prod_{i=1}^r\frac{(e_iz_iq/h_s;q)_{N_s}}{(az_iq/h_s;q)_{N_s}}\\\times
\sum_{\begin{smallmatrix}l_1,\dots,l_r\ge 0\\
0\le |{\mathbf l}|\le N_s\end{smallmatrix}}
\Bigg(\prod_{1\le i<j\le r}
\left(\frac{e_iz_iq^{l_i}-e_jz_jq^{l_j}}{e_iz_i-e_jz_j}\right)
\prod_{i,j=1}^r\frac{(e_iz_i/az_j;q)_{l_i}}
{(qe_iz_i/e_jz_j;q)_{l_i}}\\\times
\prod_{i=1}^r\frac{(ez_i;q)_{l_i}}{(e_iz_iq/h_s;q)_{l_i}}
\cdot\frac{(q^{-N_s};q)_{|{\mathbf l}|}}{(Eh_sq^{-N_s}/a^r;q)_{|{\mathbf l}|}}
\,q^{|{\mathbf l}|}\\\times
\prod_{j=1}^{s-1}\left[\frac{(Fq/h_j;q)_{N_j}}{(q/h_j;q)_{N_j}}
\prod_{i=1}^r\frac{(az_iq/f_ih_j;q)_{N_j}}{(az_iq/h_j;q)_{N_j}}\right]\\\times
\frac{(Fq/d,aq/cF;q)_{\infty}}{(aq/c,q/d;q)_{\infty}}
\prod_{i,j=1}^r\frac{(qz_i/z_j,az_iq^{1-l_j}/e_jf_iz_j,
f_jz_iq/b_iz_j;q)_{\infty}}
{(qf_jz_i/f_iz_j,qz_i/b_iz_j,az_iq^{1-l_j}/e_jz_j;q)_{\infty}}\\\times
\prod_{i=1}^r
\frac{(az_iq,q/az_i,Fq^{1-l_i}/e_iz_i,az_iq/b_iF,az_iq/df_i,
f_iq/cz_i;q)_{\infty}}
{(Ff_iq/az_i,az_iq/Ff_i,az_iq/b_i,q^{1-l_i}/e_iz_i,q/cz_i,
az_iq/d;q)_{\infty}}\\\times
\sum_{k_1,\dots,k_r=-\infty}^{\infty}
\prod_{1\le i<j\le r}\left(\frac{f_iq^{k_i}/z_i-f_jq^{k_j}/z_j}
{f_i/z_i-f_j/z_j}\right)
\prod_{i=1}^r\left(\frac{1-Ff_iq^{k_i+|{\mathbf k}|}/az_i}{1-Ff_i/az_i}\right)
\\\times
\prod_{i,j=1}^r\frac{(e_jf_iz_jq^{l_j}/az_i;q)_{k_i}}
{(f_iz_jq/b_jz_i;q)_{k_i}}
\prod_{i=1}^r\frac{(b_iF/az_i;q)_{|{\mathbf k}|}}
{(Fq^{1-l_i}/e_iz_i;q)_{|{\mathbf k}|}}\\\times
\prod_{i=1}^r\frac{(df_i/az_i,f_ih_1/az_i,\dots,f_ih_{s-1}/az_i;q)_{k_i}}
{(f_iq/cz_i,f_ih_1q^{-N_1}/az_i,\dots,f_ih_{s-1}q^{-N_{s-1}}/az_i;q)_{k_i}}\\
\times
\frac{(cF/a,Fq^{1+N_1}/h_1,\dots,Fq^{1+N_{s-1}}/h_{s-1};q)_{|{\mathbf k}|}}
{(Fq/d,Fq/h_1,\dots,Fq/h_{s-1};q)_{|{\mathbf k}|}}
\left(\frac{a^{r+1}q^{1-(N_1+\dots+N_{s-1})-|{\mathbf l}|}}
{BcdE}\right)^{|{\mathbf k}|}\Bigg)\\
=\frac{(a^rq/Eh_s;q)_{N_s}}{(q/h_s;q)_{N_s}}
\prod_{i=1}^r\frac{(e_iz_iq/h_s;q)_{N_s}}{(az_iq/h_s;q)_{N_s}}
\prod_{j=1}^{s-1}\left[\frac{(Fq/h_j;q)_{N_j}}{(q/h_j;q)_{N_j}}
\prod_{i=1}^r\frac{(az_iq/f_ih_j;q)_{N_j}}{(az_iq/h_j;q)_{N_j}}\right]\\\times
\frac{(Fq/d,aq/cF;q)_{\infty}}{(aq/c,q/d;q)_{\infty}}
\prod_{i,j=1}^r\frac{(qz_i/z_j,az_iq/e_jf_iz_j,f_jz_iq/b_iz_j;q)_{\infty}}
{(qf_jz_i/f_iz_j,qz_i/b_iz_j,az_iq/e_jz_j;q)_{\infty}}\\\times
\prod_{i=1}^r\frac{(az_iq,q/az_i,Fq/e_iz_i,az_iq/b_iF,az_iq/df_i,
f_iq/cz_i;q)_{\infty}}
{(Ff_iq/az_i,az_iq/Ff_i,az_iq/b_i,q/e_iz_i,q/cz_i,az_iq/d;q)_{\infty}}\\\times
\sum_{k_1,\dots,k_r=-\infty}^{\infty}
\Bigg(\prod_{1\le i<j\le r}\left(\frac{f_iq^{k_i}/z_i-f_jq^{k_j}/z_j}
{f_i/z_i-f_j/z_j}\right)
\prod_{i=1}^r\left(\frac{1-Ff_iq^{k_i+|{\mathbf k}|}/az_i}{1-Ff_i/az_i}\right)
\\\times
\prod_{i,j=1}^r\frac{(e_jf_iz_j/az_i;q)_{k_i}}
{(f_iz_jq/b_jz_i;q)_{k_i}}
\prod_{i=1}^r\frac{(b_iF/az_i;q)_{|{\mathbf k}|}}
{(Fq/e_iz_i;q)_{|{\mathbf k}|}}\\\times
\prod_{i=1}^r\frac{(df_i/az_i,f_ih_1/az_i,\dots,f_ih_{s-1}/az_i;q)_{k_i}}
{(f_iq/cz_i,f_ih_1q^{-N_1}/az_i,\dots,f_ih_{s-1}q^{-N_{s-1}}/az_i;q)_{k_i}}\\
\times
\frac{(cF/a,Fq^{1+N_1}/h_1,\dots,Fq^{1+N_{s-1}}/h_{s-1};q)_{|{\mathbf k}|}}
{(Fq/d,Fq/h_1,\dots,Fq/h_{s-1};q)_{|{\mathbf k}|}}
\left(\frac{a^{r+1}q^{1-(N_1+\dots+N_{s-1})}}
{BcdE}\right)^{|{\mathbf k}|}\\\times
\sum_{\begin{smallmatrix}l_1,\dots,l_r\ge 0\\
0\le |{\mathbf k}|\le N_s\end{smallmatrix}}
\prod_{1\le i<j\le r}
\left(\frac{e_iz_iq^{l_i}-e_jz_jq^{l_j}}{e_iz_i-e_jz_j}\right)
\prod_{i,j=1}^r\frac{(e_if_jz_iq^{k_j}/az_j;q)_{l_i}}
{(qe_iz_i/e_jz_j;q)_{l_i}}\\\times
\prod_{i=1}^r\frac{(ez_iq^{-|{\mathbf k}|}/F;q)_{l_i}}{(e_iz_iq/h_s;q)_{l_i}}
\cdot\frac{(q^{-N_s};q)_{|{\mathbf l}|}}{(Eh_sq^{-N_s}/a^r;q)_{|{\mathbf l}|}}
\,q^{|{\mathbf l}|}\Bigg).
\end{multline*}
Now, we evaluate the inner multiple sum by the
$a_i\mapsto e_if_iq^{k_i}/a$, $b\mapsto q^{-|{\mathbf k}|}/F$,
$c\mapsto q/h_s$, $x_i\mapsto e_iz_i$, $i=1,\dots,r$, and $N\mapsto N_s$
case of the $A_{r-1}$ $q$-Pfaff--Saalsch\"utz summation in Theorem~\ref{r32},
and obtain
\begin{multline*}
\prod_{j=1}^{s}\left[\frac{(Fq/h_j;q)_{N_j}}{(q/h_j;q)_{N_j}}
\prod_{i=1}^r\frac{(az_iq/f_ih_j;q)_{N_j}}{(az_iq/h_j;q)_{N_j}}\right]\\\times
\frac{(Fq/d,aq/cF;q)_{\infty}}{(aq/c,q/d;q)_{\infty}}
\prod_{i,j=1}^r\frac{(qz_i/z_j,az_iq/e_jf_iz_j,f_jz_iq/b_iz_j;q)_{\infty}}
{(qf_jz_i/f_iz_j,qz_i/b_iz_j,az_iq/e_jz_j;q)_{\infty}}\\\times
\prod_{i=1}^r\frac{(az_iq,q/az_i,Fq/e_iz_i,az_iq/b_iF,az_iq/df_i,
f_iq/cz_i;q)_{\infty}}
{(Ff_iq/az_i,az_iq/Ff_i,az_iq/b_i,q/e_iz_i,q/cz_i,az_iq/d;q)_{\infty}}\\\times
\sum_{k_1,\dots,k_r=-\infty}^{\infty}
\Bigg(\prod_{1\le i<j\le r}\left(\frac{f_iq^{k_i}/z_i-f_jq^{k_j}/z_j}
{f_i/z_i-f_j/z_j}\right)
\prod_{i=1}^r\left(\frac{1-Ff_iq^{k_i+|{\mathbf k}|}/az_i}{1-Ff_i/az_i}\right)
\\\times
\prod_{i,j=1}^r\frac{(e_jf_iz_j/az_i;q)_{k_i}}
{(f_iz_jq/b_jz_i;q)_{k_i}}
\prod_{i=1}^r\frac{(b_iF/az_i;q)_{|{\mathbf k}|}}
{(Fq/e_iz_i;q)_{|{\mathbf k}|}}\\\times
\prod_{i=1}^r\frac{(df_i/az_i,f_ih_1/az_i,\dots,f_ih_s/az_i;q)_{k_i}}
{(f_iq/cz_i,f_ih_1q^{-N_1}/az_i,\dots,f_ih_sq^{-N_s}/az_i;q)_{k_i}}\\
\times
\frac{(cF/a,Fq^{1+N_1}/h_1,\dots,Fq^{1+N_s}/h_s;q)_{|{\mathbf k}|}}
{(Fq/d,Fq/h_1,\dots,Fq/h_s;q)_{|{\mathbf k}|}}
\left(\frac{a^{r+1}q^{1-|N|}}
{BcdE}\right)^{|{\mathbf k}|}\Bigg),
\end{multline*}
which is the right side of \eqref{mchugl}.
\end{proof}

A special case of Theorem~\ref{mchu} immediately gives the following
summation formula as a corollary. It is an $A_{r-1}$ extension of
an identity due to W.~C.~Chu~\cite[Theorem~2]{chuwp}.

\begin{Corollary}[A multilateral very-well-poised $A_{r-1}$ $q$-IPD type
summation]\label{mchuc}
Let $a$, $b_1,\dots,b_r$, $d$, $e_1,\dots,e_r$, $z_1,\dots,z_r$, and
$h_1,\dots,h_s$ be indeterminate, let $N_1,\dots,N_s$ be nonnegative
integers, let $|N|=\sum_{i=1}^s N_i$, $r\ge 1$, and
suppose that none of the denominators in \eqref{mchucgl} vanishes. Then
\begin{multline}\label{mchucgl}
{}_{6+2s}\Psi_{6+2s}^{(r)}\bigg[a;b_1,\dots,b_r;\frac aB,d;e_1,\dots,e_r;
z_1,\dots,z_r;\\
\left.\frac{aq^{1+N_1}}{h_1},\dots,\frac{aq^{1+N_s}}{h_s};h_1,\dots,h_s
\bigg|\,q,\frac{a^rq^{1-|N|}}{dE}\right]\\
=\frac{(Bq/d,q;q)_{\infty}}{(Bq,q/d;q)_{\infty}}
\prod_{i=1}^r\frac{(az_iq,q/az_i,Bq/e_iz_i,az_iq/db_i;q)_{\infty}}
{(az_iq/b_i,q/e_iz_i,Bq/az_i,az_iq/d;q)_{\infty}}\\\times
\prod_{i,j=1}^r\frac{(qz_i/z_j,az_iq/e_jb_iz_j;q)_{\infty}}
{(qz_i/b_iz_j,az_iq/e_jz_j;q)_{\infty}}
\prod_{j=1}^s\left[\frac{(Bq/h_j;q)_{N_j}}{(q/h_j;q)_{N_j}}
\prod_{i=1}^r\frac{(az_iq/b_ih_j;q)_{N_j}}{(az_iq/h_j;q)_{N_j}}\right],
\end{multline}
provided $|a^rq^{1-|N|}/dE|<1$.
\end{Corollary}
\begin{proof}
In \eqref{mchugl}, we let $c\to a/B$ and $f_i\to b_i$, for $i=1,\dots,r$.
In this case the ${}_{6+2s}\Psi_{6+2s}^{(r)}$ series on the right side
terminates from below, and from above, and evalutes to one. In particular,
the appearance of the factor
\begin{equation*}
\prod_{i,j=1}^r(b_iz_jq/b_jz_i;q)_{k_i}^{-1}
\end{equation*}
makes the terms in the series vanish unless $k_i\ge 0$, $i=1,\dots,r$.
Similarly, the appearance of the factor
\begin{equation*}
(1;q)_{|{\mathbf k}|}
\end{equation*}
ensures that if $|{\mathbf k}|>0$, the terms of the series are zero.
In total, only the term where $k_1=\dots=k_r=0$ survives, and that term
is just one.
\end{proof}

A further specialization of Corollary~\ref{mchuc}, namely the case
$e_i\to a$, $i=1,\dots,r$, yields an $r$-dimensional generalization
of G.~Gasper's~\cite[Eq.~(5.13)]{gaseltf} very-well-poised
${}_{6+2s}\phi_{5+2s}$ summation.

\begin{Corollary}[A very-well-poised $A_{r-1}$ $q$-IPD type summation]
\label{mgasc}
Let $a$, $b_1,\dots,b_r$, $d$, $z_1,\dots,z_r$, and
$h_1,\dots,h_s$ be indeterminate, let $N_1,\dots,N_s$ be nonnegative
integers, let $|N|=\sum_{i=1}^s N_i$, $r\ge 1$, and
suppose that none of the denominators in \eqref{mgascgl} vanishes. Then
\begin{multline}\label{mgascgl}
\sum_{k_1,\dots,k_r=0}^{\infty}
\Bigg(\prod_{1\le i<j\le r}
\left(\frac{z_iq^{k_i}-z_jq^{k_j}}{z_i-z_j}\right)
\prod_{i=1}^r\left(\frac{1-az_iq^{k_i+|{\mathbf k}|}}{1-az_i}\right)\\\times
\prod_{i,j=1}^r\frac{(b_jz_i/z_j;q)_{k_i}}{(qz_i/z_j;q)_{k_i}}
\prod_{i=1}^r\frac{(az_i;q)_{|{\mathbf k}|}}{(az_iq/b_i;q)_{|{\mathbf k}|}}
\\\times
\prod_{i=1}^r\frac{(az_i/B,az_iq^{1+N_1}/h_1,\dots,az_iq^{1+N_s}/h_s;q)_{k_i}}
{(az_iq/d,az_iq/h_1,\dots,az_iq/h_s;q)_{k_i}}\\\times
\frac{(d,h_1,\dots,h_s;q)_{|{\mathbf k}|}}
{(Bq,h_1q^{-N_1},\dots,h_sq^{-N_s};q)_{|{\mathbf k}|}}\,
\left(\frac{q^{1-|N|}}d\right)^{|{\mathbf k}|}\Bigg)\\
=\frac{(Bq/d,q;q)_{\infty}}{(Bq,q/d;q)_{\infty}}
\prod_{i=1}^r\frac{(az_iq,az_iq/b_id;q)_{\infty}}
{(az_iq/b_i,az_iq/d;q)_{\infty}}\\\times
\prod_{j=1}^s\left[\frac{(Bq/h_j;q)_{N_j}}{(q/h_j;q)_{N_j}}
\prod_{i=1}^r\frac{(az_iq/b_ih_j;q)_{N_j}}{(az_iq/h_j;q)_{N_j}}\right],
\end{multline}
provided $|q^{1-|N|}/d|<1$.
\end{Corollary}

To derive the multilateral $q$-IPD type
transformation in Theorem~\ref{mchua} we need the following lemma,
which is easily established by applying Theorem~\ref{r11} twice.

\begin{Lemma}\label{11lem}
Let $a_1,\dots,a_r$, $b_1,\dots,b_r$, $c_1,\dots,c_r$, $x_1,\dots,x_r$,
$y_1,\dots,y_r$, and $z$ be indeterminate, let $r\ge 1$, and
suppose that none of the denominators in \eqref{11lemgl} vanishes. Then
\begin{multline}\label{11lemgl}
\sum_{k_1,\dots,k_r=-\infty}^{\infty}
\prod_{1\le i<j\le r}\left(\frac{x_iq^{k_i}-x_jq^{k_j}}{x_i-x_j}\right)
\prod_{i,j=1}^r\frac{(a_jx_i/x_j;q)_{k_i}}{(b_jx_i/x_j;q)_{k_i}}\,
z^{|{\mathbf k}|}\\
=\frac{(Az,q/Az;q)_{\infty}}{(Azq^r/C,Cq^{1-r}/Az;q)_{\infty}}
\prod_{i,j=1}^r\frac{(qx_i/x_j,b_jx_i/a_ix_j,
c_iy_i/a_iy_j,b_jy_iq/c_jy_j;q)_{\infty}}
{(qy_i/y_j,b_jc_iy_i/a_ic_jy_j,qx_i/a_ix_j,b_jx_i/x_j;q)_{\infty}}\\\times
\sum_{k_1,\dots,k_r=-\infty}^{\infty}
\prod_{1\le i<j\le r}\left(\frac{y_iq^{k_i}-y_jq^{k_j}}{y_i-y_j}\right)
\prod_{i,j=1}^r\frac{(a_jy_iq/c_jy_j;q)_{k_i}}{(b_jy_iq/c_jy_j;q)_{k_i}}\,
z^{|{\mathbf k}|},
\end{multline}
where $|Bq^{1-r}/A|<|z|<1$.
\end{Lemma}

We have
\begin{Theorem}[A multilateral $A_{r-1}$ $q$-IPD type transformation]
\label{mchua}
Let $a_1,\dots,a_r$, $b_1,\dots,b_r$, $c_1,\dots,c_r$, $x_1,\dots,x_r$,
$y_1,\dots,y_r$, $h_{11},\dots,h_{rs}$, and $z$ be indeterminate, let
$N_{11},\dots,N_{rs}$ be nonnegative integers, let $r\ge 1$, and
suppose that none of the denominators in \eqref{mchuagl} vanishes. Then
\begin{multline}\label{mchuagl}
\sum_{k_1,\dots,k_r=-\infty}^{\infty}
\Bigg(\prod_{1\le i<j\le r}\left(\frac{x_iq^{k_i}-x_jq^{k_j}}{x_i-x_j}\right)
\prod_{i,j=1}^r\frac{(a_jx_i/x_j;q)_{k_i}}{(b_jx_i/x_j;q)_{k_i}}\\\times
\prod_{i=1}^r\prod_{j=1}^s\frac{(h_{ij}q^{N_{ij}};q)_{|{\mathbf k}|}}
{(h_{ij};q)_{|{\mathbf k}|}}\,z^{|{\mathbf k}|}\Bigg)\\
=\prod_{i=1}^r\prod_{j=1}^s\frac{(h_{ij}q^r/C;q)_{N_{ij}}}
{(h_{ij};q)_{N_{ij}}}\cdot
\frac{(Az,q/Az;q)_{\infty}}{(Azq^r/C,Cq^{1-r}/Az;q)_{\infty}}\\\times
\prod_{i,j=1}^r\frac{(qx_i/x_j,b_jx_i/a_ix_j,
c_iy_i/a_iy_j,b_jy_iq/c_jy_j;q)_{\infty}}
{(qy_i/y_j,b_jc_iy_i/a_ic_jy_j,qx_i/a_ix_j,b_jx_i/x_j;q)_{\infty}}\\\times
\sum_{k_1,\dots,k_r=-\infty}^{\infty}
\Bigg(\prod_{1\le i<j\le r}\left(\frac{y_iq^{k_i}-y_jq^{k_j}}{y_i-y_j}\right)
\prod_{i,j=1}^r\frac{(a_jy_iq/c_jy_j;q)_{k_i}}
{(b_jy_iq/c_jy_j;q)_{k_i}}\\\times
\prod_{i=1}^r\prod_{j=1}^s\frac{(h_{ij}q^{r+N_{ij}}/C;q)_{|{\mathbf k}|}}
{(h_{ij}q^r/C;q)_{|{\mathbf k}|}}\,
z^{|{\mathbf k}|}\Bigg),
\end{multline}
provided $\left|Bq^{1-r-\sum_{i,j}N_{ij}}/A\right|<|z|<1$.
\end{Theorem}

\begin{proof}
We proceed by induction on $s$. For $s=0$ \eqref{mchuagl} is true
by Lemma~\ref{11lem}. So, suppose that the transformation is already shown
for $s \mapsto s-1$. Then (again using some elementary identities from
\cite[Appendix~I]{grhyp}),
\begin{multline*}
\sum_{k_1,\dots,k_r=-\infty}^{\infty}\Bigg(
\prod_{1\le i<j\le r}\left(\frac{x_iq^{k_i}-x_jq^{k_j}}{x_i-x_j}\right)
\prod_{i,j=1}^r\frac{(a_jx_i/x_j;q)_{k_i}}{(b_jx_i/x_j;q)_{k_i}}\\\times
\prod_{i=1}^r\prod_{j=1}^{s-1}\frac{(h_{ij}q^{N_{ij}};q)_{|{\mathbf k}|}}
{(h_{ij};q)_{|{\mathbf k}|}}\,
z^{|{\mathbf k}|}\cdot
\prod_{i=1}^r\frac{(h_{is}q^{N_{is}};q)_{|{\mathbf k}|}}
{(h_{is};q)_{|{\mathbf k}|}}\Bigg)\\
=\prod_{i=1}^r\frac 1{(h_{is};q)_{N_i}}
\sum_{k_1,\dots,k_r=-\infty}^{\infty}\Bigg(
\prod_{1\le i<j\le r}\left(\frac{x_iq^{k_i}-x_jq^{k_j}}{x_i-x_j}\right)
\prod_{i,j=1}^r\frac{(a_jx_i/x_j;q)_{k_i}}{(b_jx_i/x_j;q)_{k_i}}\\\times
\prod_{i=1}^r\prod_{j=1}^{s-1}\frac{(h_{ij}q^{N_{ij}};q)_{|{\mathbf k}|}}
{(h_{ij};q)_{|{\mathbf k}|}}\,
z^{|{\mathbf k}|}\cdot
\prod_{i=1}^r(h_{is}q^{|{\mathbf k}|};q)_{N_{is}}\Bigg).
\end{multline*}
Now we expand the last factors (those involving $(\cdot;q)_{N_{is}}$)
by applying the $x_i\mapsto h_{is}$, $n_i\mapsto N_{is}$, $i=1,\dots,r$,
and $z\mapsto q^{|{\mathbf k}|+(N_{1s}+\dots+N_{rs})}$ case of the $A_{r-1}$
summation in Theorem~\ref{r10}. We obtain
\begin{multline*}
\prod_{i=1}^r\frac 1{(h_{is};q)_{N_i}}
\sum_{k_1,\dots,k_r=-\infty}^{\infty}\Bigg(
\prod_{1\le i<j\le r}\left(\frac{x_iq^{k_i}-x_jq^{k_j}}{x_i-x_j}\right)
\prod_{i,j=1}^r\frac{(a_jx_i/x_j;q)_{k_i}}{(b_jx_i/x_j;q)_{k_i}}\\\times
\prod_{i=1}^r\prod_{j=1}^{s-1}\frac{(h_{ij}q^{N_{ij}};q)_{|{\mathbf k}|}}
{(h_{ij};q)_{|{\mathbf k}|}}\,
z^{|{\mathbf k}|}\\\times
\sum_{\begin{smallmatrix}0\le l_i\le N_{is}\\
i=1,\dots,r\end{smallmatrix}}
\prod_{1\le i<j\le r}
\left(\frac{h_{is}q^{l_i}-h_{js}q^{l_j}}{h_{is}-h_{js}}\right)
\prod_{i,j=1}^r\frac{(q^{-N_{js}}h_{is}/h_{js};q)_{l_i}}
{(qh_{is}/h_{js};q)_{l_i}}\\\times
\prod_{i=1}^rh_{is}^{l_i}\cdot
q^{|{\mathbf k}||{\mathbf l}|+(N_{1s}+\dots+N_{rs})|{\mathbf l}|
-{\binom{|{\mathbf l}|}2}+\sum_{i=1}^r{\binom{l_i}2}}\Bigg)\\
=\prod_{i=1}^r\frac 1{(h_{is};q)_{N_i}}
\sum_{\begin{smallmatrix}0\le l_i\le N_{is}\\
i=1,\dots,r\end{smallmatrix}}
\Bigg(\prod_{1\le i<j\le r}
\left(\frac{h_{is}q^{l_i}-h_{js}q^{l_j}}{h_{is}-h_{js}}\right)
\prod_{i,j=1}^r\frac{(q^{-N_{js}}h_{is}/h_{js};q)_{l_i}}
{(qh_{is}/h_{js};q)_{l_i}}\\\times
\prod_{i=1}^rh_{is}^{l_i}\cdot
q^{(N_{1s}+\dots+N_{rs})|{\mathbf l}|
-{\binom{|{\mathbf l}|}2}+\sum_{i=1}^r{\binom{l_i}2}}\\\times
\sum_{k_1,\dots,k_r=-\infty}^{\infty}
\prod_{1\le i<j\le r}\left(\frac{x_iq^{k_i}-x_jq^{k_j}}{x_i-x_j}\right)
\prod_{i,j=1}^r\frac{(a_jx_i/x_j;q)_{k_i}}{(b_jx_i/x_j;q)_{k_i}}\\\times
\prod_{i=1}^r\prod_{j=1}^{s-1}\frac{(h_{ij}q^{N_{ij}};q)_{|{\mathbf k}|}}
{(h_{ij};q)_{|{\mathbf k}|}}
\left(zq^{|{\mathbf l}|}\right)^{|{\mathbf k}|}\Bigg).
\end{multline*}
By the $z\mapsto zq^{|{\mathbf l}|}$ case of the inductive hypothesis
we obtain
\begin{multline*}
\prod_{i=1}^r\frac 1{(h_{is};q)_{N_i}}
\sum_{\begin{smallmatrix}0\le l_i\le N_{is}\\
i=1,\dots,r\end{smallmatrix}}
\Bigg(\prod_{1\le i<j\le r}
\left(\frac{h_{is}q^{l_i}-h_{js}q^{l_j}}{h_{is}-h_{js}}\right)
\prod_{i,j=1}^r\frac{(q^{-N_{js}}h_{is}/h_{js};q)_{l_i}}
{(qh_{is}/h_{js};q)_{l_i}}\\\times
\prod_{i=1}^rh_{is}^{l_i}\cdot
q^{(N_{1s}+\dots+N_{rs})|{\mathbf l}|
-{\binom{|{\mathbf l}|}2}+\sum_{i=1}^r{\binom{l_i}2}}\\\times
\prod_{i=1}^r\prod_{j=1}^{s-1}\frac{(h_{ij}q^r/C;q)_{N_{ij}}}
{(h_{ij};q)_{N_{ij}}}\cdot
\frac{(Azq^{|{\mathbf l}|},q^{1-|{\mathbf l}|}/Az;q)_{\infty}}
{(Azq^{r+|{\mathbf l}|}/C,Cq^{1-r-|{\mathbf l}|}/Az;q)_{\infty}}\\\times
\prod_{i,j=1}^r\frac{(qx_i/x_j,b_jx_i/a_ix_j,
c_iy_i/a_iy_j,b_jy_iq/c_jy_j;q)_{\infty}}
{(qy_i/y_j,b_jc_iy_i/a_ic_jy_j,qx_i/a_ix_j,b_jx_i/x_j;q)_{\infty}}\\\times
\sum_{k_1,\dots,k_r=-\infty}^{\infty}
\prod_{1\le i<j\le r}\left(\frac{y_iq^{k_i}-y_jq^{k_j}}{y_i-y_j}\right)
\prod_{i,j=1}^r\frac{(a_jy_iq/c_jy_j;q)_{k_i}}
{(b_jy_iq/c_jy_j;q)_{k_i}}\\\times
\prod_{i=1}^r\prod_{j=1}^{s-1}\frac{(h_{ij}q^{r+N_{ij}}/C;q)_{|{\mathbf k}|}}
{(h_{ij}q^r/C;q)_{|{\mathbf k}|}}
\left(zq^{|{\mathbf l}|}\right)^{|{\mathbf k}|}\Bigg)\\
=\prod_{i=1}^r\frac 1{(h_{is};q)_{N_i}}
\prod_{i=1}^r\prod_{j=1}^{s-1}\frac{(h_{ij}q^r/C;q)_{N_{ij}}}
{(h_{ij};q)_{N_{ij}}}\cdot
\frac{(Az,q/Az;q)_{\infty}}
{(Azq^r/C,Cq^{1-r}/Az;q)_{\infty}}\\\times
\prod_{i,j=1}^r\frac{(qx_i/x_j,b_jx_i/a_ix_j,
c_iy_i/a_iy_j,b_jy_iq/c_jy_j;q)_{\infty}}
{(qy_i/y_j,b_jc_iy_i/a_ic_jy_j,qx_i/a_ix_j,b_jx_i/x_j;q)_{\infty}}\\\times
\sum_{k_1,\dots,k_r=-\infty}^{\infty}\Bigg(
\prod_{1\le i<j\le r}\left(\frac{y_iq^{k_i}-y_jq^{k_j}}{y_i-y_j}\right)
\prod_{i,j=1}^r\frac{(a_jy_iq/c_jy_j;q)_{k_i}}
{(b_jy_iq/c_jy_j;q)_{k_i}}\\\times
\prod_{i=1}^r\prod_{j=1}^{s-1}\frac{(h_{ij}q^{r+N_{ij}}/C;q)_{|{\mathbf k}|}}
{(h_{ij}q^r/C;q)_{|{\mathbf k}|}}\,
z^{|{\mathbf k}|}\\\times
\sum_{\begin{smallmatrix}0\le l_i\le N_{is}\\
i=1,\dots,r\end{smallmatrix}}
\prod_{1\le i<j\le r}
\left(\frac{h_{is}q^{l_i}-h_{js}q^{l_j}}{h_{is}-h_{js}}\right)
\prod_{i,j=1}^r\frac{(q^{-N_{js}}h_{is}/h_{js};q)_{l_i}}
{(qh_{is}/h_{js};q)_{l_i}}\\\times
\prod_{i=1}^rh_{is}^{l_i}\cdot \left(\frac {q^r}C\right)^{|{\mathbf l}|}
q^{|{\mathbf k}||{\mathbf l}|+(N_{1s}+\dots+N_{rs})|{\mathbf l}|
-{\binom{|{\mathbf l}|}2}+\sum_{i=1}^r{\binom{l_i}2}}\Bigg).
\end{multline*}

Now, we evaluate the inner multiple sum by the $x_i\mapsto h_{is}$,
$n_i\mapsto N_{is}$, $i=1,\dots,r$,
and $z\mapsto q^{r+|{\mathbf k}|+(N_{1s}+\dots+N_{rs})}/C$ case of the $A_{r-1}$
summation in Theorem~\ref{r10}. We obtain
\begin{multline*}
\prod_{i=1}^r\frac 1{(h_{is};q)_{N_i}}
\prod_{i=1}^r\prod_{j=1}^{s-1}\frac{(h_{ij}q^r/C;q)_{N_{ij}}}
{(h_{ij};q)_{N_{ij}}}\cdot
\frac{(Az,q/Az;q)_{\infty}}
{(Azq^r/C,Cq^{1-r}/Az;q)_{\infty}}\\\times
\prod_{i,j=1}^r\frac{(qx_i/x_j,b_jx_i/a_ix_j,
c_iy_i/a_iy_j,b_jy_iq/c_jy_j;q)_{\infty}}
{(qy_i/y_j,b_jc_iy_i/a_ic_jy_j,qx_i/a_ix_j,b_jx_i/x_j;q)_{\infty}}\\\times
\sum_{k_1,\dots,k_r=-\infty}^{\infty}\Bigg(
\prod_{1\le i<j\le r}\left(\frac{y_iq^{k_i}-y_jq^{k_j}}{y_i-y_j}\right)
\prod_{i,j=1}^r\frac{(a_jy_iq/c_jy_j;q)_{k_i}}
{(b_jy_iq/c_jy_j;q)_{k_i}}\\\times
\prod_{i=1}^r\prod_{j=1}^{s-1}\frac{(h_{ij}q^{r+N_{ij}}/C;q)_{|{\mathbf k}|}}
{(h_{ij}q^r/C;q)_{|{\mathbf k}|}}\,
z^{|{\mathbf k}|}
\cdot\prod_{i=1}^r(h_iq^{r+|{\mathbf k}|}/C;q)_{N_{is}}\Bigg),
\end{multline*}
which, after an elementary manipulation of $q$-shifted factorials,
gives us the right side of \eqref{mchuagl}, as desired.
\end{proof}

Finally, we provide four more multilateral transformations of
$q$-IPD type. Unfortunately, we were not able
to find multiple extensions of Propositions~\ref{km1} or \ref{km0}
which are as deep as the identities in Theorems~\ref{mchu} and \ref{mchua}.
The following theorems are obtained by combining Propositions~\ref{km1}
and \ref{km0} each with Lemmas~\ref{lem1} and \ref{lem2}, thus giving rise
to four different multilateral transformations.

\begin{Theorem}[A multilateral $A_{r-1}$ $q$-IPD type transformation]
\label{mkm1}
Let $a$, $b$, $c_1,\dots,c_r$, $d$, $e_1,\dots,e_r$, $x_1,\dots,x_r$,
$y_1,\dots,y_r$, and $h_1,\dots,h_s$ be indeterminate, let $N$ be an integer,
let $m_1,\dots,m_s$ be nonnegative integers, let $|m|=\sum_{i=1}^s m_i$,
$r\ge 1$, and
suppose that none of the denominators in \eqref{mkm1gl} vanishes. Then
\begin{multline}\label{mkm1gl}
\sum_{k_1,\dots,k_r=-\infty}^{\infty}
\Bigg(\prod_{1\le i<j\le r}\left(\frac {x_iq^{k_i}-x_jq^{k_j}}{x_i-x_j}\right)
\prod_{i,j=1}^r(c_jx_i/x_j;q)_{k_i}^{-1}
\prod_{i=1}^rx_i^{rk_i-|{\mathbf k}|}\\\times
(-1)^{(r-1)|{\mathbf k}|}
q^{-\binom{|{\mathbf k}|}2+r\sum_{i=1}^r\binom{k_i}2}\\\times
\frac{(a,b,h_1q^{m_1},\dots,h_sq^{m_s};q)_{|{\mathbf k}|}}
{(d,h_1,\dots,h_s;q)_{|{\mathbf k}|}}
\left(\frac{Eq^{1-r-N}}{ab}\right)^{|{\mathbf k}|}\Bigg)\\
=(Eq^{-r})^N\frac{(Eq^{1-r}/a,Eq^{1-r}/b,dq^r/E;q)_{\infty}}
{(q/a,q/b,d;q)_{\infty}}\\\times
\prod_{i=1}^s\frac{(h_iq^r/E;q)_{m_i}}{(h_i;q)_{m_i}}
\prod_{i,j=1}^r\frac{(qx_i/x_j,c_jy_iq/e_jy_j;q)_{\infty}}
{(qy_i/y_j,c_jx_i/x_j;q)_{\infty}}\\\times
\sum_{k_1,\dots,k_r=-\infty}^{\infty}
\Bigg(\prod_{1\le i<j\le r}\left(\frac {y_iq^{k_i}-y_jq^{k_j}}{y_i-y_j}\right)
\prod_{i,j=1}^r(c_jy_iq/e_jy_j;q)_{k_i}^{-1}
\prod_{i=1}^ry_i^{rk_i-|{\mathbf k}|}\\\times
(-1)^{(r-1)|{\mathbf k}|}
q^{-\binom{|{\mathbf k}|}2+r\sum_{i=1}^r\binom{k_i}2}\\\times
\frac{(aq^r/E,bq^r/E,h_1q^{r+m_1}/E,\dots,h_sq^{r+m_s}/E;q)_{|{\mathbf k}|}}
{(dq^r/E,h_1q^r/E,\dots,h_sq^r/E;q)_{|{\mathbf k}|}}
\left(\frac{Eq^{1-r-N}}{ab}\right)^{|{\mathbf k}|}\Bigg),
\end{multline}
provided $|Eq^{1-r}/ab|<|q^N|<|Eq^{|m|}/Cd|$.
\end{Theorem}
\begin{proof}
We have, for
$|Eq^{1-r}/ab|<|q^N|<|Eq^{|m|}/Cd|$,
\begin{multline}\label{22t1gl}
{}_{2+s}\psi_{2+s}\!\left[\begin{matrix}a,b,h_1q^{m_1},\dots,h_sq^{m_s}\\
Cq^{1-r},d,h_1,\dots,h_s\end{matrix}\,;q,\frac{Eq^{1-r-N}}{ab}\right]\\
=(Eq^{-r})^N\frac{(Eq^{1-r}/a,Eq^{1-r}/b,Cq/E,dq^r/E;q)_{\infty}}
{(q/a,q/b,Cq^{1-r},d;q)_{\infty}}
\prod_{i=1}^s\frac{(h_iq^r/E;q)_{m_i}}{(h_i;q)_{m_i}}\\\times
{}_{2+s}\psi_{2+s}\!\left[\begin{matrix}aq^r/E,bq^r/E,h_1q^{r+m_1}/E,\dots,
h_sq^{r+m_s}/E\\
Cq/E,dq^r/E,h_1q^r/E,\dots,h_sq^r/E\end{matrix}
\,;q,\frac{Eq^{1-r-N}}{ab}\right],
\end{multline}
by the $q$-IPD type transformation in \eqref{km1gl}.
Now we apply Lemma~\ref{lem1} to the $_{2+s}\psi_{2+s}$'s on the left and
on the right side of this transformation.
Specifically, we rewrite the $_{2+s}\psi_{2+s}$ on left side of \eqref{22t1gl}
by the $b_i\mapsto c_i$, $i=1,\dots,r$, and
\begin{equation*}
f(n)=\frac{(a,b,h_1q^{m_1},\dots,h_sq^{m_s};q)_n}{(d,h_1,\dots,h_s;q)_n}
\left(\frac{Eq^{1-r-N}}{ab}\right)^n
\end{equation*}
case of Lemma~\ref{lem1}.
The $_{2+s}\psi_{2+s}$ on the right side of \eqref{22t1gl} is rewritten
by the $b_i\mapsto c_iq/e_i$, $x_i\mapsto y_i$, $i=1,\dots,r$, and
\begin{equation*}
f(n)=\frac{(aq^r/E,bq^r/E,h_1q^{r+m_1}/E,\dots,h_sq^{r+m_s}/E;q)_n}
{(dq^r/E,h_1q^r/E,\dots,h_sq^r/E;q)_n}
\left(\frac{Eq^{1-r-N}}{ab}\right)^n
\end{equation*}
case of Lemma~\ref{lem1}.
Finally, we divide both sides of the resulting equation by
\begin{equation}\label{divbs}
\frac{(q;q)_{\infty}}{(Cq^{1-r};q)_{\infty}}
\prod_{i,j=1}^r\frac{(c_jx_i/x_j;q)_{\infty}}{(qx_i/x_j;q)_{\infty}}
\end{equation}
and simplify to obtain \eqref{mkm1gl}.
\end{proof}

\begin{Theorem}[A multilateral $A_{r-1}$ $q$-IPD type transformation]
\label{mkm1a}
Let $a_1,\dots,a_r$, $b$, $c_1,\dots,c_r$, $d$, $e_1,\dots,e_r$,
$x_1,\dots,x_r$, $y_1,\dots,y_r$, and $h_1,\dots,h_s$ be indeterminate,
let $N$ be an integer, let $m_1,\dots,m_s$ be nonnegative integers,
let $|m|=\sum_{i=1}^s m_i$, $r\ge 1$, and suppose that none of the
denominators in \eqref{mkm1agl} vanishes. Then
\begin{multline}\label{mkm1agl}
\sum_{k_1,\dots,k_r=-\infty}^{\infty}
\Bigg(\prod_{1\le i<j\le r}\left(\frac {x_iq^{k_i}-x_jq^{k_j}}{x_i-x_j}\right)
\prod_{i,j=1}^r\frac{(a_jx_i/x_j;q)_{k_i}}{(c_jx_i/x_j;q)_{k_i}}\\\times
\frac{(b,h_1q^{m_1},\dots,h_sq^{m_s};q)_{|{\mathbf k}|}}
{(d,h_1,\dots,h_s;q)_{|{\mathbf k}|}}
\left(\frac{Eq^{1-r-N}}{Ab}\right)^{|{\mathbf k}|}\Bigg)\\
=(Eq^{-r})^N\frac{(Eq^{1-r}/b,dq^r/E;q)_{\infty}}
{(q/b,d;q)_{\infty}}
\prod_{i=1}^s\frac{(h_iq^r/E;q)_{m_i}}{(h_i;q)_{m_i}}\\\times
\prod_{i,j=1}^r
\frac{(qx_i/x_j,c_jx_i/a_ix_j,c_jy_iq/e_jy_j,e_iy_i/a_iy_j;q)_{\infty}}
{(qy_i/y_j,c_je_iy_i/a_ie_jy_j,c_jx_i/x_j,x_iq/a_ix_j;q)_{\infty}}\\\times
\sum_{k_1,\dots,k_r=-\infty}^{\infty}
\Bigg(\prod_{1\le i<j\le r}\left(\frac {y_iq^{k_i}-y_jq^{k_j}}{y_i-y_j}\right)
\prod_{i,j=1}^r\frac{(a_jy_iq/e_jy_j;q)_{k_i}}{(c_jy_iq/e_jy_j;q)_{k_i}}
\\\times
\frac{(bq^r/E,h_1q^{r+m_1}/E,\dots,h_sq^{r+m_s}/E;q)_{|{\mathbf k}|}}
{(dq^r/E,h_1q^r/E,\dots,h_sq^r/E;q)_{|{\mathbf k}|}}
\left(\frac{Eq^{1-r-N}}{Ab}\right)^{|{\mathbf k}|}\Bigg),
\end{multline}
provided $|Eq^{1-r}/Ab|<|q^N|<|Eq^{|m|}/Cd|$.
\end{Theorem}
\begin{proof}
We have, for
$|Eq^{1-r}/Ab|<|q^N|<|Eq^{|m|}/Cd|$,
\begin{multline}\label{22t2gl}
{}_{2+s}\psi_{2+s}\!\left[\begin{matrix}A,b,h_1q^{m_1},\dots,h_sq^{m_s}\\
Cq^{1-r},d,h_1,\dots,h_s\end{matrix}\,;q,\frac{Eq^{1-r-N}}{Ab}\right]\\
=(Eq^{-r})^N\frac{(Eq^{1-r}/A,Eq^{1-r}/b,Cq/E,dq^r/E;q)_{\infty}}
{(q/A,q/b,Cq^{1-r},d;q)_{\infty}}
\prod_{i=1}^s\frac{(h_iq^r/E;q)_{m_i}}{(h_i;q)_{m_i}}\\\times
{}_{2+s}\psi_{2+s}\!\left[\begin{matrix}Aq^r/E,bq^r/E,h_1q^{r+m_1}/E,\dots,
h_sq^{r+m_s}/E\\
Cq/E,dq^r/E,h_1q^r/E,\dots,h_sq^r/E\end{matrix}
\,;q,\frac{Eq^{1-r-N}}{Ab}\right],
\end{multline}
by the $q$-IPD type transformation in \eqref{km1gl}.
Now we apply Lemma~\ref{lem2} to the $_{2+s}\psi_{2+s}$'s on the left and
on the right side of this transformation.
Specifically, we rewrite the $_{2+s}\psi_{2+s}$ on left side of \eqref{22t2gl}
by the $b_i\mapsto c_i$, $i=1,\dots,r$, and
\begin{equation*}
g(n)=\frac{(b,h_1q^{m_1},\dots,h_sq^{m_s};q)_n}{(d,h_1,\dots,h_s;q)_n}
\left(\frac{Eq^{1-r-N}}{Ab}\right)^n
\end{equation*}
case of Lemma~\ref{lem2}.
The $_{2+s}\psi_{2+s}$ on the right side of \eqref{22t2gl} is rewritten
by the $a_i\mapsto a_iq/e_i$, $b_i\mapsto c_iq/e_i$, $x_i\mapsto y_i$,
$i=1,\dots,r$, and
\begin{equation*}
g(n)=\frac{(bq^r/E,h_1q^{r+m_1}/E,\dots,h_sq^{r+m_s}/E;q)_n}
{(dq^r/E,h_1q^r/E,\dots,h_sq^r/E;q)_n}
\left(\frac{Eq^{1-r-N}}{Ab}\right)^n
\end{equation*}
case of Lemma~\ref{lem2}.
Finally, we divide both sides of the resulting equation by
\begin{equation}\label{divbs2}
\frac{(q,Cq^{1-r}/A;q)_{\infty}}{(Cq^{1-r},q/A;q)_{\infty}}
\prod_{i,j=1}^r\frac{(c_jx_i/x_j,x_iq/a_ix_j;q)_{\infty}}
{(qx_i/x_j,c_jx_i/a_ix_j;q)_{\infty}}
\end{equation}
and simplify to obtain \eqref{mkm1agl}.
\end{proof}

\begin{Theorem}[A multilateral $A_{r-1}$ $q$-IPD type transformation]
\label{mkm0}
Let $a_1,\dots,a_r$, $b$, $c_1,\dots,c_r$, $d$, $e_1,\dots,e_r$,
$x_1,\dots,x_r$, $y_1,\dots,y_r$, and $h_1,\dots,h_s$ be indeterminate,
let $N$ be an integer, let $m_1,\dots,m_s$ be nonnegative integers,
let $|m|=\sum_{i=1}^s m_i$, $r\ge 1$, and
suppose that none of the denominators in \eqref{mkm0gl} vanishes. Then
\begin{multline}\label{mkm0gl}
\sum_{k_1,\dots,k_r=-\infty}^{\infty}
\Bigg(\prod_{1\le i<j\le r}\left(\frac {x_iq^{k_i}-x_jq^{k_j}}{x_i-x_j}\right)
\prod_{i,j=1}^r(c_jx_i/x_j;q)_{k_i}^{-1}
\prod_{i=1}^rx_i^{rk_i-|{\mathbf k}|}\\\times
(-1)^{(r-1)|{\mathbf k}|}
q^{-\binom{|{\mathbf k}|}2+r\sum_{i=1}^r\binom{k_i}2}\\\times
\frac{(A,b,h_1q^{m_1},\dots,h_sq^{m_s};q)_{|{\mathbf k}|}}
{(d,h_1,\dots,h_s;q)_{|{\mathbf k}|}}
\left(\frac{Eq^{1-r-N}}{Ab}\right)^{|{\mathbf k}|}\Bigg)\\
=(Eq^{-r})^N\frac{(Eq^{1-r}/b,Cq/E,dq^r/E;q)_{\infty}}
{(q/A,q/b,d;q)_{\infty}}\\\times
\prod_{i=1}^s\frac{(h_iq^r/E;q)_{m_i}}{(h_i;q)_{m_i}}
\prod_{i,j=1}^r\frac{(qx_i/x_j,e_jy_i/a_jy_j;q)_{\infty}}
{(qy_i/y_j,c_jx_i/x_j;q)_{\infty}}\\\times
\sum_{k_1,\dots,k_r=-\infty}^{\infty}
\Bigg(\prod_{1\le i<j\le r}\left(\frac {y_iq^{k_i}-y_jq^{k_j}}{y_i-y_j}\right)
\prod_{i,j=1}^r(e_jy_i/a_jy_j;q)_{k_i}^{-1}
\prod_{i=1}^ry_i^{rk_i-|{\mathbf k}|}\\\times
(-1)^{(r-1)|{\mathbf k}|}
q^{-\binom{|{\mathbf k}|}2+r\sum_{i=1}^r\binom{k_i}2}\\\times
\frac{(E/C,Eq^{1-r}/d,Eq^{1-r}/h_1,\dots,Eq^{1-r}/h_s;q)_{|{\mathbf k}|}}
{(Eq^{1-r}/b,Eq^{1-r-m_1}/h_1,\dots,Eq^{1-r-m_s}/h_s;q)_{|{\mathbf k}|}}
\left(\frac{Cdq^{N-|m|}}E\right)^{|{\mathbf k}|}\Bigg),
\end{multline}
provided $|Eq^{1-r}/Ab|<|q^N|<|Eq^{|m|}/Cd|$.
\end{Theorem}
\begin{proof}
We have, for
$|Eq^{1-r}/Ab|<|q^N|<|Eq^{|m|}/Cd|$,
\begin{multline}\label{22t3gl}
{}_{2+s}\psi_{2+s}\!\left[\begin{matrix}A,b,h_1q^{m_1},\dots,h_sq^{m_s}\\
Cq^{1-r},d,h_1,\dots,h_s\end{matrix}\,;q,\frac{Eq^{1-r-N}}{Ab}\right]\\
=(Eq^{-r})^N\frac{(Eq^{1-r}/A,Eq^{1-r}/b,Cq/E,dq^r/E;q)_{\infty}}
{(q/A,q/b,Cq^{1-r},d;q)_{\infty}}
\prod_{i=1}^s\frac{(h_iq^r/E;q)_{m_i}}{(h_i;q)_{m_i}}\\\times
{}_{2+s}\psi_{2+s}\!\left[\begin{matrix}E/C,Eq^{1-r}/d,Eq^{1-r}/h_1,\dots,
Eq^{1-r}/h_s\\
Eq^{1-r}/A,Eq^{1-r}/b,Eq^{1-r-m_1}/h_1,\dots,Eq^{1-r-m_s}/h_s\end{matrix}
\,;q,\frac{Cdq^{N-|m|}}E\right],
\end{multline}
by the $q$-IPD type transformation in \eqref{km0gl}.
Now we apply Lemma~\ref{lem1} to the $_{2+s}\psi_{2+s}$'s on the left and
on the right side of this transformation.
Specifically, we rewrite the $_{2+s}\psi_{2+s}$ on left side of \eqref{22t3gl}
by the $b_i\mapsto c_i$, $i=1,\dots,r$, and
\begin{equation*}
f(n)=\frac{(A,b,h_1q^{m_1},\dots,h_sq^{m_s};q)_n}{(d,h_1,\dots,h_s;q)_n}
\left(\frac{Eq^{1-r-N}}{Ab}\right)^n
\end{equation*}
case of Lemma~\ref{lem1}.
The $_{2+s}\psi_{2+s}$ on the right side of \eqref{22t3gl} is rewritten
by the $b_i\mapsto c_i/a_i$, $x_i\mapsto y_i$, $i=1,\dots,r$, and
\begin{equation*}
f(n)=\frac{(E/C,Eq^{1-r}/d,Eq^{1-r}/h_1,\dots,Eq^{1-r}/h_s;q)_n}
{(Eq^{1-r}/b,Eq^{1-r-m_1}/h_1,\dots,Eq^{1-r-m_s}/h_s;q)_n}
\left(\frac{Cdq^{N-|m|}}E\right)^n
\end{equation*}
case of Lemma~\ref{lem1}.
Finally, we divide both sides of the resulting equation by \eqref{divbs}
and simplify to obtain \eqref{mkm0gl}.
\end{proof}

\begin{Theorem}[A multilateral $A_{r-1}$ $q$-IPD type transformation]
\label{mkm0a}
Let $a_1,\dots,$ $a_r$, $b$, $c_1,\dots,c_r$, $d$, $e_1,\dots,e_r$,
$x_1,\dots,x_r$, $y_1,\dots,y_r$, and $h_1,\dots,h_s$ be indeterminate,
let $N$ be an integer, let $m_1,\dots,m_s$ be nonnegative integers,
let $|m|=\sum_{i=1}^s m_i$, $r\ge 1$, and suppose that none of the
denominators in \eqref{mkm0agl} vanishes. Then
\begin{multline}\label{mkm0agl}
\sum_{k_1,\dots,k_r=-\infty}^{\infty}
\Bigg(\prod_{1\le i<j\le r}\left(\frac {x_iq^{k_i}-x_jq^{k_j}}{x_i-x_j}\right)
\prod_{i,j=1}^r\frac{(a_jx_i/x_j;q)_{k_i}}{(c_jx_i/x_j;q)_{k_i}}\\\times
\frac{(b,h_1q^{m_1},\dots,h_sq^{m_s};q)_{|{\mathbf k}|}}
{(d,h_1,\dots,h_s;q)_{|{\mathbf k}|}}
\left(\frac{Eq^{1-r-N}}{Ab}\right)^{|{\mathbf k}|}\Bigg)\\
=(Eq^{-r})^N\frac{(Eq^{1-r}/b,dq^r/E;q)_{\infty}}
{(q/b,d;q)_{\infty}}
\prod_{i=1}^s\frac{(h_iq^r/E;q)_{m_i}}{(h_i;q)_{m_i}}\\\times
\prod_{i,j=1}^r
\frac{(qx_i/x_j,c_jx_i/a_ix_j,e_jy_i/a_jy_j,c_iy_iq/e_iy_j;q)_{\infty}}
{(qy_i/y_j,c_ie_jy_i/a_je_iy_j,c_jx_i/x_j,x_iq/a_ix_j;q)_{\infty}}\\\times
\sum_{k_1,\dots,k_r=-\infty}^{\infty}
\Bigg(\prod_{1\le i<j\le r}\left(\frac {y_iq^{k_i}-y_jq^{k_j}}{y_i-y_j}\right)
\prod_{i,j=1}^r\frac{(e_jy_i/c_jy_j;q)_{k_i}}{(e_jy_i/a_jy_j;q)_{k_i}}
\\\times
\frac{(Eq^{1-r}/d,Eq^{1-r}/h_1,\dots,Eq^{1-r}/h_s;q)_{|{\mathbf k}|}}
{(Eq^{1-r}/b,Eq^{1-r-m_1}/h_1,\dots,Eq^{1-r-m_s}/h_s;q)_{|{\mathbf k}|}}
\left(\frac{Cdq^{N-|m|}}E\right)^{|{\mathbf k}|}\Bigg),
\end{multline}
provided $|Eq^{1-r}/Ab|<|q^N|<|Eq^{|m|}/Cd|$.
\end{Theorem}
\begin{proof}
We have, for
$|Eq^{1-r}/Ab|<|q^N|<|Eq^{|m|}/Cd|$, \eqref{22t3gl}
by the $q$-IPD type transformation in \eqref{km0gl}.
Now we apply Lemma~\ref{lem2} to the $_{2+s}\psi_{2+s}$'s on the left and
on the right side of this transformation.
Specifically, we rewrite the $_{2+s}\psi_{2+s}$ on left side of \eqref{22t3gl}
by the $b_i\mapsto c_i$, $i=1,\dots,r$, and
\begin{equation*}
g(n)=\frac{(b,h_1q^{m_1},\dots,h_sq^{m_s};q)_n}{(d,h_1,\dots,h_s;q)_n}
\left(\frac{Eq^{1-r-N}}{Ab}\right)^n
\end{equation*}
case of Lemma~\ref{lem2}.
The $_{2+s}\psi_{2+s}$ on the right side of \eqref{22t3gl} is rewritten
by the $a_i\mapsto e_i/c_i$, $b_i\mapsto e_i/a_i$, $x_i\mapsto y_i$,
$i=1,\dots,r$, and
\begin{equation*}
g(n)=\frac{(Eq^{1-r}/d,Eq^{1-r}/h_1,\dots,Eq^{1-r}/h_s;q)_n}
{(Eq^{1-r}/b,Eq^{1-r-m_1}/h_1,\dots,Eq^{1-r-m_s}/h_s;q)_n}
\left(\frac{Cdq^{N-|m|}}E\right)^n
\end{equation*}
case of Lemma~\ref{lem2}.
Finally, we divide both sides of the resulting equation by \eqref{divbs2}
and simplify to obtain \eqref{mkm0agl}.
\end{proof}

The $e_i=a_iq$, $i=1,\dots,r$, cases of Theorems~\ref{mkm0} and \ref{mkm0a}
yield two $A_{r-1}$ extensions of W.~C.~Chu's~\cite[Eq.~15]{chubs} bilateral
transformation. If we specialize these identities further by setting
$c_i=q$, $i=1,\dots,r$, we obtain two $A_{r-1}$ extensions of
G.~Gasper's~\cite[Eq.~(19)]{gassum} $q$-IPD type transformation.

\section*{Acknowledgements}
We wish to thank Professors George Gasper and Mourad Ismail for
valuable suggestions which led to an improvement of the article.

\end{document}